\documentclass{article}
\usepackage{graphicx} 
\usepackage{amsmath}
\usepackage{amssymb}
\usepackage[margin=30truemm]{geometry}
\usepackage{amsthm}
\usepackage{color}
\usepackage{enumerate}
\theoremstyle{definition}
\newtheorem{thm}{Theorem}[section]
\newtheorem{prop}[thm]{Proposition}
\newtheorem{lemma}[thm]{Lemma}
\newtheorem{definition}[thm]{Definition}
\newtheorem{rem}[thm]{Remark}
\newtheorem*{pf}{Proof}

\title{Higher-dimensional multifractal analysis for the cusp winding process  on hyperbolic surfaces}
\author{Yuya Arima\\
 Graduate School of Mathematics, Nagoya University,\\
Furocho, Chikusaku, Nagoya, 464-8602, Japan\\
(e-mail: yuya.arima.c0@math.nagoya-u.ac.jp)}

\begin{document}

\maketitle

\begin{abstract}
    We perform a multifractal analysis of the growth rate of the number of cusp windings for the geodesic flow on hyperbolic surfaces with $m \geq 1$ cusps. Our main theorem establishes a conditional variational principle for the Hausdorff dimension spectrum of the multi-cusp winding process. Moreover, we show that  the dimension spectrum defined on $\mathbb{R}_{>0}^m$ is real analytic. To prove the main theorem we use a countable Markov shift with a finitely primitive transition matrix and thermodynamic formalism.   
\end{abstract}
\section{Introduction}
Throughout this paper, we consider the Poincar\'e disc model $(\mathbb{D}, d)$ of two-dimensional hyperbolic space. For details on hyperbolic geometry, we refer the reader to \cite{Dal}, \cite{katok1992fuchsian} and \cite{beardon2012geometry}.
Denote by Conf$(\mathbb{D})$ the set of orientation-preserving isometries of $(\mathbb{D},d)$. A subgroup $G$ of $\text{Conf}(\mathbb{D})$ is called a Fuchsian group if $G$ is a discrete subgroup of the group Conf$(\mathbb{D})$. Fuchsian groups play an important role in the uniformization of hyperbolic surfaces and geometric group
theory. (see \cite[Section 6]{beardon2012geometry}).

We always assume that $G$ is a  generalized Schottky group with $m\ge1$ parabolic generators. We refer to Section 2.1 for the precise definition. Note that $G$ is a non-elementary finitely generated free Fuchsian group with respect to a set of generators $G_0$ (see \cite[Proposition 1.6]{Dal}).  The set $G_0$ is given by $G_0=\Gamma_0\cup H_0$ where $\Gamma_0:=\{\gamma_1^{\pm 1},\cdots,\gamma_m^{\pm 1}\}$ is a set of parabolic generators,  and $H_0:=\{h_1^{\pm1},\cdots,h_n^{\pm1}\}$ a set of $n\ge 1$ hyperbolic generators. 
 We define the limit set of $G$ as $\Lambda(G):=\overline{\bigcup_{g\in G}\{g0\}}\setminus\bigcup_{g\in G}\{g0\} \subset \partial{\mathbb{D}}$. A limit point $x\in\Lambda(G)$ is called a conical limit point if there exists $(g_n)_{n\in\mathbb{N}}\subset G$ such that $\lim_{n\to\infty}g_n0=x$ and the sequence $(\inf_{z\in[0,x)}d(g_n0,z))_{n\in\mathbb{N}}$ is bounded where $[0,x)$ denote the geodesic ray connecting $0$ and $x$. By \cite[\text{Theorem}\ 10.2.5]  {beardon2012geometry}, we can completely decompose the limit set into the  set of conical limit points and the set of non-conical limit points. More precisely, the conical limit set $\Lambda_c(G)$ is given by 
\begin{align}\label{eq decomposition of limit set}
    \Lambda_c(G)=\Lambda(G)\setminus\bigcup_{g\in G }\bigcup_{i=1}^m\{gp_i\},
\end{align}
where $p_i$ is the unique fixed point of $\gamma_i$ ($i=1,2,\cdots,m$).
Hence, $\Lambda(G)\setminus\Lambda_c(G)$ is a countable set.
A conical limit point $x\in \Lambda_c(G)$ is coded by infinite sequences over the set {$G_0$} as follows. 
Let $s_x$ be the oriented geodesic ray from 0 towards $x\in\Lambda_c(G)$ and $R\subset \mathbb{D}$ denote the Dirichlet fundamental domain for $G$ centered at $0$. The oriented geodesic $s_x$ intersects infinitely many copies $R$, $g_0(x)R$, $g_0(x)g_1(x)R$,... of $R$, with $g_i(x)\in G_0$ and $i\in\mathbb{N}$. We obtain the infinite sequence $\omega(x)=g_0(x)g_1(x)g_2(x)\cdots\in G_0^{\mathbb{N}}$, which is necessarily reduced,  that is $g_{i-1}(x)g_i(x)\neq I$ for all $i\in\mathbb{N}$. 

We recall the definition of the multi-cusp winding process from \cite{Jaerisch2016AMA} (see also \cite{jaerisch2021mixed}). For a conical limit point $x\in\Lambda_c(G)$ the sequence $\omega(x)=g_0(x)g_1(x)g_2(x)\cdots$ defines a block sequence $B_i(x)$,  $i\in \mathbb{N}$, such that
\begin{align*}
\omega(x)=B_1(x)B_2(x)\cdots    
\end{align*}
 where each $B_i(x)$ is either a hyperbolic generator, or a maximal block of consecutive appearances of the same parabolic generator.
 By construction, for some $\gamma\in\Gamma_0$, $l\geq0$ and $i\in\mathbb{N}$ a block $B_i(x)=\gamma^{l+1}$ means that the projection of $s_x$ onto $\mathbb{D}/G$ winding $l$ times around the cusp associated with $\gamma$. Motivated by this, we define the multi-cusp winding process
$(a_{i,j})_{1\leq i\leq m, j\geq1}:\Lambda_c(G)\rightarrow \mathbb{N}\cup\{0\}$ as 
\begin{align*}
    a_{i,j}(x)= \left\{
 \begin{array}{cc}
   l   & \text{if}\  B_j(x)=\gamma_{i}^{l+1}\ \text{or}\ B_j(x)=\gamma_{i}^{-(l+1)},\ l\geq1\\
   0   & \text{otherwise}
 \end{array}
 \right..
\end{align*}
For $\boldsymbol{\alpha}=(\alpha_1,\alpha_2,\cdots,\alpha_m)\in[0,\infty]^{m}$ we define the level set by 
\begin{align*}
    J(\boldsymbol{\alpha}):=\left\{x\in\Lambda_c(G):\frac{1}{n}\sum_{j=0}^{n-1}a_{i,j}(x)=\alpha_i,\ 1\leq i\leq m\right\}.
\end{align*}
We have a multifractal decomposition of the conical limit set 
\begin{align*}
    \Lambda_c(G)=\left(\bigcup_{\boldsymbol{\alpha}\in[0,\infty]^m}J(\boldsymbol{\alpha}) \right)\cup J_{\text{ir}}
\end{align*}
where $J_{\text{ir}}$ is the irregular set defined by 
\begin{align*}
    J_{\text{ir}}:=\left\{x\in\Lambda_c(G):\exists i\in\{1,2,\cdots m\}\ \text{s.t.\ }\lim_{n\to\infty}\frac{1}{n}\sum_{j=0}^{n-1}a_{i,j}\ \text{does\ not\ exist}\right\}.
\end{align*}
By \cite[Proposition 2.1]{jaerisch2021mixed}, for $\boldsymbol{\alpha}\in[0,\infty]^m$ the set $J(\boldsymbol{\alpha})$ is a dense subset of  $\Lambda_c(G)$. Thus, this decomposition of the $\Lambda_c(G)$ is quite complicated. To investigate the fine structure of $\Lambda_c(G)$, for each $\boldsymbol{\alpha}\in[0,\infty]^m$ it is necessary to analyze $J(\boldsymbol{\alpha})$. To do this, we define the multi-cusp winding spectrum 
\begin{align*}
    b:[0,\infty]^m\rightarrow [0,\dim_H\Lambda(G)],\ b(\boldsymbol{\alpha})=\dim_HJ(\boldsymbol{\alpha})
\end{align*}
where $\dim_H(\cdot)$ denotes the Hausdorff dimension. Also, the function $b$ is called simply the dimension spectrum. The dimension spectrum describes the growth rate of the number of cusp windings. For all $\boldsymbol{\alpha}\in [0,\infty]^m$ the box dimension of $J(\boldsymbol{\alpha})$ is the box dimension of $\Lambda_c(G)$ since $J(\boldsymbol{\alpha})$ is a dense subset of $\Lambda_c(G)$. This is the reason for using Hausdorff dimension instead of box dimension to measure the size of $J(\boldsymbol{\alpha})$.

Denote by $f$ the Bowen-Series map associated with $G$ and $G_0$ and by $\tilde f$ the induced system derived from the Bowen-Series map $f$ (see Section 2.1 for the precise definition). By (1), the maximal $\tilde f$-invariant set is the conical limit set $\Lambda_c(G)$.  Denote the set of $\tilde f|_{\Lambda_c(G)}$-invariant Borel probability  measures by $M(\tilde f)$. For $\mu\in M(\tilde f)$ we denote the measure-theoretic entropy as $h(\mu)$ and we define the Lyapunov exponent of $\mu\in M(\tilde f)$ as $\lambda(\mu):=\int\log|\tilde f'|d\mu $. We define the Hausdorff dimension of a $\tilde f|_{\Lambda_c(G)}$-invariant Borel probability  measure $\mu\in M(\tilde f)$ as
\begin{align*}
    \dim_H(\mu):=\inf\{\dim_H(Z):\mu(Z)=1\}.
\end{align*}

\begin{thm}\label{main}
     Let $G$ be a generalized Schottky group with $m\geq1$ parabolic generators. Then the following holds:
    
    (1) For $\boldsymbol{\alpha}=(\alpha_1,\cdots,\alpha_m)\in(0,\infty)^m$ we have 
    \begin{align*}
       \ \  b(\boldsymbol{\alpha})=\max\left\{\frac{h(\mu)}{\lambda(\mu)}:\mu\in M(\tilde f),\ \lambda(\mu)<\infty,\ \int a_{i,1}d\mu=\alpha_i,\ 1\leq i\leq m\right\}.
    \end{align*}
    
    (2) There exists $\mu\in M(\tilde f)$ such that 
    \begin{align*}
    \dim_H(\mu)=\dim_H J(\infty,\infty,\cdots,\infty)=\dim_H(\Lambda(G)).
\end{align*}

    (3) The dimension spectrum $b$ is real-analytic on $(0,\infty)^m$.\\
    
    (4) We have $\dim_H(J_{\text{ir}})=\dim_H(\Lambda(G))$.\\
\end{thm}
If $m=1$ then we can say a little more.
\begin{prop}\label{Proposition m=1 case}
Let $G$ be a generalized Schottky group with one parabolic generator. Then the dimension spectrum $b$ is strictly increasing on $(0,\infty)$ and we have 
    \begin{align*}
    \lim_{\alpha\to\infty}\dim_H(J(\alpha))=1.
       \end{align*}    
\end{prop}
    
We use the same setting as the section 2 in \cite{jaerisch2021mixed}. Our main claim is item $(3)$ in Theorem $\ref{main}$. The induced map $\tilde f$ is coded by a countable Markov shift with a finitely primitive transition matrix. Moreover, $\tilde f$ is a uniformly expanding map which satisfies the Renyi condition and the bounded distortion property. The multi-cusp winding process defines a weakly H\"older potential. Thus, Theorem \ref{main} is similar to the 
Iommi and Jordan's results on the conditional variational principle and regularity of a Birkhoff spectrum \cite[Theorem 1.3, and Theorem 1.4]{iommi2015multifractal}. The biggest difference between Theorem \ref{main} and \cite[Theorem 1.3, and Theorem 1.4]{iommi2015multifractal} is that the dimension spectrum is defined on the $m$-dimensional Euclidean space $\mathbb{R}^m$. Moreover, the Markov maps  considered in \cite{iommi2015multifractal} are coded by a full-shift on a countable alphabet, but  $\tilde f$ can not be coded by a full-shift. Also, the potential considered in \cite{iommi2015multifractal} does not take the value zero,  while the multi-cusp winding process can take zero values. Therefore, even in the case of $m=1$, we cannot directly apply the \cite[Theorem 1.3, and Theorem 1.4]{iommi2015multifractal}. Barreira, Saussol, and Schmeling \cite{barreira2002higher} perform a higher-dimensional version of multifractal analysis for a hyperbolic dynamical system and they show the conditional variational principle and regularity of a Birkhoff spectrum defined on higher-dimensional Euclidean space for a weakly H\"older potential. However, to our knowledge, there is no known example and result on 
the regularity of a dimension spectrum defined on higher-dimensional Euclidean space for a Markov map with countable branches and an unbounded weakly H\"older potential.  
Also, item $(1)$ of Theorem \ref{main} extends the conditional variational formula for the cusp winding spectrum \cite[Proposition 2.1]{jaerisch2021mixed}. We refer to  \cite{barany2021birkhoff,barreira2002higher,barreira2000sets,climenhaga2012topological,fan2015multifractal,iommi2015multifractal,iommi2015multifractale,iommi2017transience,jaerisch2021mixed,johansson2008multifractal,Jordan2017BirkhoffSF,Kessebhmer2007HigherdimensionalMV,rush2023multifractal} for 
other results of the thermodynamic formalism related to Theorem \ref{main}.

In the case of $m=1$, there is a study closely related to Theorem \ref{main}.
 Using the Gauss map $T:(0,1)\setminus\mathbb{Q}\rightarrow(0,1)\setminus\mathbb{Q}$ defined by $T(x)=1/x$ mod $1$, Iommi and Jordan \cite[Corollary 6.6]{iommi2015multifractal} showed that the conditional variational principle for the Birkhoff spectrum of the arithmetic mean of the regular continued fraction digits, the Birkhoff spectrum is real-analytic on $(1,\infty)$. Moreover, the Birkhoff spectrum is strictly increasing on $(1,\infty)$ and tends to $1$ as a variable goes to infinity. 
 We can associate geodesics on the quotient of the hyperbolic surface by the modular group with the regular continued fraction \cite{series1985modular}. Thus, Theorem \ref{main} can be regarded as an analog of \cite[Corollary 6.6]{iommi2015multifractal} when $m=1$. In the case $m\geq 1$, the cusp winding spectrum was also studied in \cite{Jaerisch2016AMA} and \cite{Munday2011OnHD}.\\

\textbf{Methods of proofs.}
The proof of Theorem \ref{main} mainly follows the arguments of Iommi and Jordan \cite{iommi2015multifractal} (see also Remark \ref{remark-gap}  below). By \cite[Proposition 2.1]{jaerisch2021mixed},  we have for all $\boldsymbol{\alpha}\in(0,\infty)^m$, 
\begin{align}\label{eq:conditional variational principle with epsilon}
 b(\boldsymbol{\alpha}) = \lim_{\epsilon \to 0}\sup \left\{ \frac{h(\mu)}{\lambda(\mu)}: \mu \in M(\tilde f),\lambda(\mu)<\infty, \left|\int a_i d \mu - \alpha_i  \right|<\epsilon,\ 1\leq i\leq m 
   \right\}.
\end{align}
We improve  (\ref{eq:conditional variational principle with epsilon}) to the conditional variational principle for the dimension spectrum, that is 
\begin{align}\label{eq conditional variational principle}
    b(\boldsymbol{\alpha})=\sup\left\{\frac{h(\mu)}{\lambda(\mu)}:\mu\in M(\tilde f),\ \lambda(\mu)<\infty,\ \int a_{i,1}d\mu=\alpha_i,\ 1\leq i\leq m\right\}
\end{align}
for $\boldsymbol{\alpha}\in(0,\infty)^m$. To do this, 
We will show that the function 
\begin{align*}
\boldsymbol{\alpha}\mapsto\sup\left\{\frac{h(\mu)}{\lambda(\mu)}:\mu\in M(\tilde f),\ \lambda(\mu)<\infty,\ \int a_{i,1}d\mu=\alpha_i,\ 1\leq i\leq m\right\}
\end{align*}
is continuous on $(0,\infty)^m$. 

\begin{rem} \label{remark-gap}
   It seems to us that there is a gap in the proof of \cite[Lemma 3.2]{iommi2015multifractal}. Namely, the existence of compactly supported invariant measures approximating  entropy, Lyapunov exponent and a given potential in the proof of \cite[Lemma 3.2]{iommi2015multifractal} requires further explanation. This gap can be fixed by  using \cite[Lemma 2.3]{jaerisch2021mixed} and \cite[Main theorem]{takahasi2020entropy}, see also Propositions \ref{conditional variational principle with epsilon} and \ref{improve}.  
\end{rem}
Next, we show that the dimension spectrum is real-analytic on $(0,\infty)^m$. Since the induced map $\tilde f$ is a Markov map with countable many branches, we can use the thermodynamic formalism in the symbolic setting. In particular, we can consider the topological pressure $P$ defined by Mauldin and Urba\'nski (see Section 2 or \cite{mauldin2003graph} for the precise definition) which has good properties. Moreover, those good properties can be translated into properties of the topological pressure with respect to the dynamical system $(\tilde f|_{\Lambda_c(G)},\Lambda_c(G))$ (see Section 2 for the precise definition).
To prove that the dimension spectrum is real-analytic, define the function $p$ by 
\begin{align*}
 p:\mathbb{R}^m\times\mathbb{R}^m\times\mathbb{R}\rightarrow \mathbb{R},\ \ \    
p(\boldsymbol{\alpha},\boldsymbol{q},b):=P(\langle\boldsymbol{q},(\Phi+\boldsymbol{\alpha})\rangle-b\log|\tilde f'|)\  
\end{align*} 
where $\langle\cdot ,\cdot \rangle$ denotes the standard inner product on $\mathbb{R}^m$ and $\Phi(x):=(-a_{1,1}(x),-a_{2,1}(x),\cdots,-a_{m,1}(x)) \nonumber$, and the function $p$ is related to the dimension spectrum as follows:
by using the above good properties for the topological pressure $P$, considering the asymptotic behavior of the function $p$, and using the improved conditional  variational formula for the dimension spectrum, 
 one can show that for $\boldsymbol{\alpha}\in(0,\infty)^m$ there exists $\boldsymbol{q}(\boldsymbol{\alpha})\in(0,\infty)^m$ such that 
\begin{align}\label{equation G}
    p(\boldsymbol{\alpha},\boldsymbol{q}(\boldsymbol{\alpha}),b(\boldsymbol{\alpha}))=0\ \text{and}\  \frac{\partial}{\partial q_i}p(\boldsymbol{\alpha},\boldsymbol{q}(\boldsymbol{\alpha}),b(\boldsymbol{\alpha}))=0\ \text{for all}\ 1\leq i \leq m.
\end{align}
Moreover, we show in Proposition \ref{prop more important proposition} that the supremum in (\ref{eq conditional variational principle}) is attained. Hence, by implicit function theorem, we can show that $b$ is real-analytic. \\

\textbf{Plan of the paper.}
In Section 2.1, we introduce the precise definition of the Bowen-Series map $f$ and the multi-cusp winding process. Moreover, we describe the induced dynamical system derived from the Bowen-Series map $\tilde f$ and the dimension spectrum. In Section 2.2, we explain the relationship between a Markov map and a countable Markov shift and then provide basic notions of the thermodynamic formalism which are often used in this paper. 

Section 3 is devoted to improving the conditional variational formula for the dimension spectrum $(\ref{eq:conditional variational principle with epsilon})$ to the conditional variational principle for the dimension spectrum $(\ref{eq conditional variational principle})$.

In Section \ref{Convergence of the topological pressure and Bowen's formula}, we give the necessary and sufficient conditions for convergence of pressure for a potential mainly used in this paper. By this condition,  Bowen's formula is shown and we obtain item $(2)$ of Theorem \ref{main}. Moreover, we prove that the maximal measure is the equilibrium measure for the geometric potential $-\dim_H(\Lambda_c(G))\log|\tilde f'|$. 

In Section \ref{Relationship between the topological pressure and dimension spectrum}, we associate the dimension spectrum with the topological pressure. Also, in the proof of Proposition \ref{important lemma}, we show that for $\boldsymbol{\alpha}\in(0,\infty)^m$ the measure which attains the supremum of the conditional variational principle for the dimension spectrum is the equilibrium measure for the potential $\langle\boldsymbol{q}(\boldsymbol{\alpha}),(\Phi+\boldsymbol{\alpha})\rangle-b(\boldsymbol{\alpha})\log|\tilde f'|$ where $\boldsymbol{q}(\boldsymbol{\alpha})$ is the value in (\ref{equation G}). 


In Section \ref{Regularity of the dimension spectrum on a open "good" set}, using the results of Section \ref{Relationship between the topological pressure and dimension spectrum}, 
 we prove that the dimension spectrum is real-analytic on $(0,\infty)^m$. Also, we prove item (4) of Theorem \ref{main}.

In Section \ref{section proof 2},  we prove Proposition \ref{Proposition m=1 case}. It is important to consider the shape of the graph of the dimension spectrum. By item (2) of Theorem \ref{main}, we know that $\dim_H(J(\infty))$ is equal to the Hausdorff dimension of the equilibrium measure for the geometric potential $-\dim_H(\Lambda_c(G))\log|\tilde f'|$ and to the Hausdorff dimension of the conical limit set $\Lambda_c(G)$. Thus, item (4) of Theorem \ref{main} states that the graph of the dimension spectrum 
$b$ is typical.

\section{Preliminaies}
\subsection{The Bowen-Series map and multi-cusp winding process}
In this section, we will first give some definitions of hyperbolic geometry and the notations used throughout this paper. 
Recall that the elements of Conf$(\mathbb{D})$ can be  classified as hyperbolic, parabolic, or elliptic using their fixed points.   A parabolic element $\gamma\in\text{Conf}(\mathbb{D})$ has exactly one fixed point $p$ which is in the Euclidean boundary $\partial \mathbb{D}$ of $\mathbb{D}$. Moreover, we have $|\gamma'(p)|=1$ where $|\cdot|$ denotes the Euclidean metric norm on $\mathbb{R}^2$.
Let $m\geq 1$ and let
$\Gamma_0:=\lbrace \gamma_1^{\pm 1},\cdots,\gamma_m^{\pm 1}\rbrace \subset \text{Conf}(\mathbb{D})$
be parabolic generators with a fixed point $p_i\ (i=1,2,\cdots,m)$ of the parabolic generator $\gamma_i^{\pm1}$. Let $n\geq1$ and let
$H_0:=\lbrace h_1^{\pm 1},\cdots,h_{n}^{\pm 1} \rbrace \subset \text{Conf}(\mathbb{D})$
be hyperbolic generators. Denote $G_0:=\Gamma_0 \cup H_0$ and write $G_0:=\lbrace g_1,g_2,\cdots,g_{2(m+n)} \rbrace$.  For $g \in \text{Conf}(\mathbb{D})$, let $\Delta(g):=\lbrace z \in \partial \mathbb{D} : |g'(z)|\geq1 \rbrace$ and put
\begin{align}
    \Delta:=\bigcup_{i=1}^{2(m+n)}\Delta({g_i}).\nonumber
\end{align}
We assume that 
    $\overline{\Delta({g_i})\cup \Delta({g_i^{-1}})}\cap \overline{{\Delta({g_j})\cup \Delta({g_j^{-1}})}}= \emptyset$
for $i,j\in\{1,2,\cdots,2(m+n)\}$ and $g_i,g_j\in G_0$ with $ g_i \neq g_j^{\pm 1}$. Then, the subgroup $G$ of Conf$(\mathbb{D})$ generated by $G_0$ is called a generalized Schottky group with $m$ parabolic generators. As mentioned in the introduction, $G$ is a non-elementary finitely generated free Fuchsian group. The surface $\mathbb{D} / G$ has $m$ cusps given by $p_i \ (i=1,2,\cdots,m)$.
\begin{definition}
    The Bowen-Series map \cite{bowen1979markov} is given by 
\begin{align}
   f:\Delta \rightarrow \partial \mathbb{D}, \ f|_{\Delta({g_i})}=g_i\ \ (i=1,\cdots,2(m+n)).\nonumber
\end{align}
\end{definition}
  For all $x\in \Lambda_c(G)$ there uniquely exists  $\omega(x)=\omega_1\omega_2\cdots \in G_0^\mathbb{N}$ such that $f^n(x)\in\Delta(\omega_n)$ and if $\omega_n\ (n\in \mathbb{N})$ is parabolic then $\exists m>n$ such that $\omega_m\neq\omega_n$. Then, $\omega$ defines a sequence of blocks $B_i$ such that
$\omega(x)=B_1(x)B_2(x)\cdots$
 where each $B_i(x)$ is either a hyperbolic generator, or a maximal block of consecutive appearances of the same parabolic generator.

\begin{definition}
    The cusp winding process $(a_{i,j})_{1\leq i\leq m,j\geq1}:\Lambda_c(G)\rightarrow \mathbb{N}_{\geq0}$ is given by
\begin{align} 
 a_{i,j}(x)= \left\{
 \begin{array}{cc}
   m   & \text{if}\  B_j(x)=\gamma_i^{\pm{(m+1)}} (\gamma_i \in \Gamma_0 )\\
   0   & \text{otherwise}\nonumber
 \end{array}
 \right. .
\end{align}
\end{definition}
Define 
\begin{align}
\mathcal{A}:=\bigcup_{l=0}^{\infty} \lbrace \gamma^{l+1}g:\gamma \in \Gamma_0,g \in H_0\rbrace \cup H_0\nonumber
\end{align}
and for $\omega=\omega_0\cdots\omega_{n-1} \in \mathcal{A}$, define the set
$\Delta({\omega}):=\Delta({\omega_0}) \cap f^{-1}\Delta({\omega_1}) \cap \cdots \cap f^{-(n-1)}\Delta({\omega_{n-1}}).$

\begin{definition}
Define the inducing time $\tau:\mathcal{A}\rightarrow \mathbb{N}$ by $\tau(\gamma^{l+1}h)={l+1}$ ($\gamma \in \Gamma_0,\ h\in H_0)$ and $\tau|_{H_0}=1$. The induced Markov map with the Markov partition $\lbrace \Delta({\omega}) \rbrace_{\omega \in A}$ is given by
\begin{align}
\tilde f: \bigcup_{\omega \in A} \Delta(\omega)  \rightarrow \partial \mathbb{D}, \ \ \tilde f|_{\Delta({\omega})}=f^{\tau(\omega)}. \nonumber
\end{align}
\end{definition}

Note that  $a_{i,j}=a_{i,1}\circ \tilde f^{j-1} $ for $j\geq 1$ and $i\in\lbrace1,2,\cdots,m\rbrace$, the maximal $\tilde f$-invariant set is the conical limit set $\Lambda_c(G)$ by (\ref{eq decomposition of limit set}), and $\tilde f|_{\Lambda_c(G)}$ is uniformly expanding map that is there exists $c>1$ such that $\inf_{x\in \Lambda_c(G)} |\tilde f'(x)|>c$ by definition of the Dirichlet fundamental domain. 
By definition, there exists a constant $Z\geq1$ such that for all $h\in H_0$ and $x\in\Delta(h)$ we have
\begin{align*}
    \frac{1}{Z}\leq |\tilde f'(x)|\leq Z.
\end{align*}
In this paper, we denote simply $\tilde f|_{\Lambda_c(G)}$ as $\tilde f$ and consider the pair $(\tilde f, \Lambda_c(G))$ as a dynamical system.

   We denote by $M({\tilde f})$ the set of $\tilde f$-invariant 
Borel probability measures on $\Lambda_c(G)$. For $\mu\in M(\tilde f)$ we denote by $h(\mu)$ the measure-theoretic entropy and by $\lambda(\mu):=\int \log|\tilde f'| d \mu$ the Lyapunov exponent. Since $\tilde f$ is uniformly expanding, for all $\mu\in M(\tilde f)$ we have $\lambda(\mu)>0$.

\begin{definition}
    Define for $\boldsymbol{\alpha}=(\alpha_1, \alpha_2,\cdots,\alpha_m) \in [0,\infty]^m$ the level sets,
    \begin{align}
    J(\boldsymbol{\alpha}):=\left\{ x \in \Lambda_c(G):\lim_{n\to \infty}\frac{1}{n}\sum_{i=0}^{n-1}a_{i,1}(\tilde f^{i}(x)) = {\alpha_i},i\in\lbrace1,2,\cdots,m\rbrace \right\},\nonumber
    \end{align}
    and the dimension spectrum 
\begin{align}
   b:[0,\infty]^m\rightarrow\mathbb{R},\ \ b(\boldsymbol{\alpha}):=\dim_HJ(\boldsymbol{\alpha}).\nonumber
\end{align}
 The dimension spectrum describes the growth rate of the number of cusp windings.
\end{definition}

\subsection{Thermodynamic formalism}
In this section, we describe the thermodynamic formalism that is often used to show the regularity of the dimension spectrum. The thermodynamic formalism usually plays an important role in multifractal analysis. For more detail information regarding the thermodynamic formalism we refer the reader to \cite[Section 2]{mauldin2003graph} and \cite[section 17]{urubanskinoninvertible}.

Recall that $\tilde f $ is a Markov map. Thus, $\tilde f $ determines a $\mathcal{A}\times \mathcal{A}$ matrix $A$ by $A_{ab}=1$ if $\Delta_b\subset\tilde f \Delta_a$ and $A_{ab}=0$ otherwise. Define
\begin{align}
    \Sigma_A:=\lbrace \omega \in \mathcal{A}^{\mathbb{N}}:A_{\omega_{n-1},\omega_{n}}=1,\ n\in \mathbb{N} \rbrace\nonumber.
\end{align}
An admissible word of length $n$ is a string $(\omega_0,\omega_1,\cdots,\omega_{n-1})\in\mathcal{A}^n$ such that $A_{\omega_{i-1},\omega_{i}}=1$ for all $i=0,\cdots,n-1$. We denote by $E^n$ the set of all admissible words of length $n$ and put $E^*:=\bigcup_{n=1}^{\infty}E^n$. For $F\subset\mathcal{A}$ and $n\in\mathbb{N}$, we define $F^n:=\lbrace\omega\in E^n:\omega_i\in F,\  0\leq i\leq n-1\rbrace$ and $F^\infty:=\lbrace\omega\in \Sigma_A:\omega_{i-1}\in F,\ i\in\mathbb{N}\rbrace$. For convenience, put $E^0=\lbrace \emptyset \rbrace$. If $\omega \in E^n$ then the corresponding cylinder set in $\Sigma_A$ is defined by $[\omega]:=\lbrace \tau \in \Sigma_A:\tau_i=\omega_i, 0\leq i\leq n-1\rbrace$. Note that  $\Sigma_A$ is finitely primitive that is there exists $n\in\mathbb{N}$ and a finite set $\Lambda\subset E^n$ such that for any $\omega, \tau\in E^*$ there exists $\lambda\in\Lambda$ such that $\omega\lambda\tau\in E^*$. This means that $\Sigma_A$ is finitely irreducible that is there exists a finite set $\Lambda\subset E^*$ such that for any $\omega, \tau\in E^*$ there exists $\lambda\in\Lambda$ such that $\omega\lambda\tau\in E^*$. \\
We endow $\Sigma_A$ with the topology generated by the cylinders. We use the metric $d:\Sigma_A\times\Sigma_A\rightarrow \mathbb{R}$ on $\Sigma_A$ defined by 
\begin{align}\label{eq: metric}
     d(\omega,\omega')= \left\{
 \begin{array}{cc}
   {2^{-k}}   & \text{if}\  \omega_i=\omega_i'\ \text{for\ all}\ i=0,\cdots,k-1\ \text{and}\ \omega_{k}\neq\omega_{k}'\\
   0   & \text{if}\ \omega=\omega'
 \end{array}
 \right. .
\end{align}
This generates the same topology as that of the cylinder sets. Recall that $\Sigma_A$ is compact with respect to the topology generated by $d$ if and only if a set of symbol is a finite set. Thus, in our case, $\Sigma_A$ is not compact.

The following lemma is used to relate the dynamical system $(\tilde f, \Lambda_c(G))$ and $(\Sigma_A,\sigma)$ where $\sigma$ denote the shift map defined on $\Sigma_A$.
\begin{prop}[{\cite[Proposition 3.1]{jaerisch2021mixed}}] 
For any $\omega=(\omega_0,\omega_1,\cdots) \in \Sigma_A$ the set 
    $\bigcap_{j=0}^{\infty}{\tilde{f}^{-j}\Delta_{\omega_j}}\nonumber$
    is a singleton.
\end{prop}
By this proposition, the coding map $\pi:\Sigma_A\rightarrow \pi(\Sigma_A)$ given by 
\begin{align}
    \pi(\omega)\in\bigcap_{n=0}^\infty{\tilde{f}^{-j}\Delta_{\omega_j}},\ \ \ (\omega=(\omega_0,\omega_1,\cdots) \in \Sigma_A)\nonumber
\end{align}
is well-defined, and bijection. Moreover, $\pi$ and $\pi^{-1}$ are continuous, and $\pi$ satisfies $\tilde{f}(\pi(\omega))=\pi(\sigma(\omega))$ for $\omega \in \Sigma_A$. Note that we have $\pi(\Sigma_A)=\Lambda_c(G)$. 
\begin{definition}\label{definition variation}
The $n$-variation of a function $\phi:\Sigma_A\rightarrow \mathbb{R}$ is defined by 
\begin{align}
    V_n(\phi):=\sup \lbrace |\phi(\omega)-\phi(\tau)|:\omega,\tau\in\Sigma_A,\ \omega_i=\tau_i\ \text{for}\ 0\leq i\leq n-1\rbrace\nonumber.
\end{align}
We say that $\phi$ has summable variation if $\sum_{n=1}^{\infty} V_n(\phi)<\infty$. Also $\phi$ is called weakly H\"older if there exists $A>0$ and $W\in(0,1)$ such that we have $V_n(\phi)\leq A W^n$.

\end{definition}
Note that $a_{i,1}\circ \pi\ (i=1,2,\cdots, m)$ and $\log|\tilde{f}'\circ \pi|$ are weakly H\"older. All potentials appear on this paper are weakly H\"older.

 \begin{definition}[{\cite{mauldin2003graph}}]\label{Definition of pressure}
Let $\phi:\Sigma_A\rightarrow \mathbb{R}$ is a continuous function. For $F\subset\mathcal{A}$ the topological pressure of $\phi$ is defined by 
 \begin{align}
     P_F(\phi):=\lim_{n\to\infty}\frac{1}{n}\log \sum_{\omega\in E^n\cap F^n}\exp\left(\sup_{\tau\in[\omega\cap F]}\sum_{j=0}^{n-1}\phi(\sigma^i(\tau))\right)\nonumber.
 \end{align}
 where $[\omega\cup F]:=\lbrace\tau\in F^\infty:\tau\in[\omega]\rbrace$. If $F=\mathcal{A}$, we simply write $[\omega]$ for $[\omega\cap F]$ and $P(\phi)$ for $P_F(\phi)$.
\end{definition}
The above limit always exits, but it can be positive infinity.

In the remainder of this section we assume that $\phi:\Sigma_A\rightarrow\mathbb{R}$ is weakly H\"older. 

The following lemma is very useful when determining the finiteness of topological pressure.
\begin{lemma}[{\cite[Proposition 2.1.9]{mauldin2003graph}}]\label{original convergence lemma}
$P(\phi)<\infty$ if and only if $\sum_{\omega\in E}\exp(\sup\phi|_{[\omega]})<\infty$. 
\end{lemma}
An important property of the topological pressure is that it can be approximated by compact sets.
\begin{prop}[{\cite[Theorem 2.1.5]{mauldin2003graph}}]\label{approximation property}
    We have
        $P(\phi)=\sup\lbrace P_F(\phi)\rbrace\nonumber$
    where the supremum is taken over all finite subsets $F$ of $\mathcal{A}$.
\end{prop}
The topological pressure satisfies the following variational principle.
\begin{thm}[{\cite[Theorem 2.1.8]{mauldin2003graph}}]
We have 
    $P(\phi)=\sup \left\{h(\mu)+\int \phi d\mu \right\},\nonumber$
    where the supremum is taken over all $\sigma$-invariant ergodic Borel probability measures $\mu$ such that $\int \phi d\mu>-\infty$ and $h(\mu)$ is the measure-theoretic entropy with respect to $\sigma$.
\end{thm}
We denote by $M(\sigma)$ the set of $\sigma$-invariant Borel probability measures on $\Sigma_A$. 
\begin{definition}
    A measure $\mu \in {M}(\sigma)$ is called an equilibrium measure for $\phi$ if 
         $P(\phi)= h(\mu)+\int \phi d\mu\nonumber.$
\end{definition}
A Gibbs measure is a measure that is related to topological pressure as follows.
\begin{definition}
    A probability measure $\mu$ is called a Gibbs measure for the potential $\phi$ if there exists a constant $M>1$ such that, for any cylinder $[\omega]\ (\omega \in E^n,\ n\in\mathbb{N})$ and any $\tau\in [\omega]$ we have 
    \begin{align}
        \frac{1}{M}\leq\frac{\mu([\omega])}{\exp(-nP(\phi)+\sum_{j=0}^{n-1}\phi(\sigma(\tau)))}\leq M\nonumber.
    \end{align}
\end{definition}
Furthermore, the topological pressure satisfies the following basic and important results.
\begin{thm}[{\cite[Theorem 2.2.9]{mauldin2003graph}}]
 Suppose that the potential $\phi$ satisfies $P(\phi)<\infty$ and\\
     $\sum_{\omega\in E}\inf(-\phi|_{[\omega]})\exp(\inf\phi|_{[\omega]})<\infty\nonumber$
 (i.e. $\int -\phi d\mu<\infty$ for all Gibbs measure $\mu\in {M}(\sigma)$). Then, there exists a unique Gibbs measure $\mu_\phi \in\mathcal{M}_\sigma$ such that 
      $P(\phi)= h(\mu_{\phi})+\int \phi d\mu_\phi.$
 \end{thm}
 The following theorem is heavy used to prove the regularity of the dimension spectrum.
 \begin{thm}
     Let $I\subset \mathbb{R}$ be an open interval satisfying the following conditions: \\
     $(1)$ $\left\{ \psi_t \right\}_{t\in I}$ is a family of weakly H\"{o}lder functions on $\Sigma_A$. \\
     $(2)$ The function $t\rightarrow \psi_t$ is analytic on $I$.\\
     $(3)$ $P(\psi_t)<\infty$ for all $t\in I$.\\
     Then $t \rightarrow P(\psi_t)$ is real analytic on $I$.
 \end{thm}
Next, we describe the thermodynamic formalism with respect to the dynamical system $(\tilde f,\Lambda_c(G))$.
\begin{definition}
The topological pressure of a continuous function $\phi : \Lambda_c(G)\rightarrow \mathbb{R}$ is defined by
\begin{equation}
P_{\tilde f}(\phi)=\sup \left\{ h(\mu) + \int \phi d\mu :  \int \phi d\mu > -\infty , \mu \in {M}({\tilde f}) \right\}.\nonumber
\end{equation}\
\end{definition}
Since $\pi$ is a bijection and $\tilde f \circ \pi = \pi \circ \sigma$, there exists a bijection between the space of $\sigma$-invariant Borel probability measures ${M}({\sigma})$ and the space of $\tilde f$-invariant Borel probability measures ${M}({\tilde f})$. Thus, we obtain $P_{\tilde f}(\phi)=P(\phi \circ \pi)$.  We will denote both pressures by $P$.

We can relate the Hausdorff dimension of the conical limit set $\Lambda_c(G)$ and the topological pressure $P$ by the following well-known formula, which is called Bowen's formula. 
\begin{thm}[{\cite[Theorem 4.2.13]{mauldin2003graph}}]\label{Bowen} We have 
    $\dim_H(\Lambda_c(G))=\inf\{t\in\mathbb{R}:P(-t\log|\tilde f'|)\leq0\}.$
\end{thm}
The Bowen's formula is improved later (see Theorem \ref{useful Bowen}).

\section{Improved conditional variational principle for the dimension spectrum}
This section is devoted to improving the conditional variational formula for the dimension spectrum (\ref{eq:conditional variational principle with epsilon}) (see Proposition \ref{conditional variational principle with epsilon}) to the conditional variational principle for the dimension spectrum (\ref{eq conditional variational principle}) in Proposition \ref{improve}.

Recall that $G$ is a non-elementary finite generated free Fuchsian group with $m\geq 1$ parabolic generators. Further, $\tilde f$ denotes the induced dynamical system derived from the Bowen-Series map $f$. We denote by $(a_{i,j})_{1\leq i\leq m,j\geq1}$ the multi-cusp winding process. For simplicity, we will use the notations $a_{i,1}=a_i$ for $i\in\lbrace1,2,\cdots,m\rbrace$. Define 
\begin{align}
    \Psi: \lbrace \mu \in {M}({\tilde f}): \lambda(\mu)<\infty \rbrace \rightarrow \mathbb{R}^m,\quad\Psi(\mu):=\left(\int a_1 d\mu,\int a_2 d\mu,\cdots,\int a_m d\mu\right).\nonumber
\end{align}

We put $\mathbb{R}^m_{>0}:=(0,\infty)^m$ and use $\text{Int}$ to denote the interior of a subset of  $\mathbb{R}^m$. 
By considering convex combinations in ${M}({\tilde f})$, we see that $\text{Int}(\text{Im}\Psi)$ is convex set.
In fact, we have the following. 
\begin{lemma}\label{imformation of U}
    We have $\text{Int}(\text{Im}\Psi)=\mathbb{R}^m_{>0}$.
\end{lemma}
\begin{pf}
 Since for all $1\leq i\leq m$ and $x\in\Lambda_c(G)$ we have $a_i(x)\geq0$, we obtain $\text{Int}(\text{Im}\Psi)\subset\mathbb{R}^m_{>0}$. Therefore, it is suffices to show that $\mathbb{R}^m_{>0}\subset \text{Int}(\text{Im}\Psi)$. For $1\leq i \leq m$ and $l\geq1$ we put $x({i,l})=\pi(\gamma_i^lh_1h_1\gamma_i^lh_1h_1\cdots)$ and $x_0=\pi(h_1h_1\cdots)$. For all $M>0$ and $1\leq i\leq m$ we have $\int a_i d\delta_{(x({i,l})+\tilde f(x({i,l})
))/2}=l/2$ and for $k\neq i$ we have $\int a_k d\delta_{(x({i,l})+\tilde f(x({i,l})))/2}=0$ where $\delta_{(x({i,l})+\tilde f(x({i,l})))/2}$ denote the point mass measure at ${(x({i,l})+\tilde f(x({i,l})))/2}$. Let $\boldsymbol{y}\in\mathbb{R}^m_{>0}$. We can take a large enough natural number $l\in\mathbb{N}$ such that  
\begin{align*}
\boldsymbol{y}\in\text{Conv}\{\Psi(\delta_{(x({i,l})+\tilde f(x({i,l})))/2}):1\leq i\leq m\}
\end{align*}
where $\text{Conv}\{\Psi(\delta_{(x({i,l})+\tilde f(x({i,l})))/2}):1\leq i\leq m\}$ denotes the convex hull of $\{\Psi(\delta_{(x({i,l})+\tilde f(x({i,l})))/2}):1\leq i\leq l\}$. Hence, by considering convex combination $\delta_{x_0}$ with $\tilde f$-invariant measures $\delta_{(x({1,l})+\tilde f(x({1,l})))/2}$, $\delta_{(x({2,l})+\tilde f(x({2,l})))/2}$, $\cdots$,$\delta_{(x({m,l})+\tilde f(x({m,l})))/2}$, we obtain $\boldsymbol{y} \in \text{Int}(\text{Im}\Psi)$ and $\mathbb{R}_{>0}^m\subset \text{Int}(\text{Im}\Psi)$.   $\qed$
\end{pf}

Next, we introduce the result of Jaerisch and Takahasi. They show the following variational formula. Recall that one of the main purpose of this paper is to improve the following variational formula and to prove the regularity of the dimension spectrum.  

\begin{prop}[Conditional variational principle with $\epsilon$,
{\cite[Proposition 2.1]{jaerisch2021mixed}}]\label{conditional variational principle with epsilon}
For $\boldsymbol{\alpha}\in[0,\infty]^m$ we have
\begin{align}
b(\boldsymbol{\alpha}) = \lim_{\epsilon \to 0}\sup \left\{ \frac{h(\mu)}{\lambda(\mu)}: \mu \in M(\tilde f),\lambda(\mu)<\infty, \left|\int a_i d \mu - \alpha_i  \right|<\epsilon,\ 1\leq i\leq m 
   \right\}.\nonumber
\end{align}
\end{prop}

To remove the $\epsilon$ in above formula, we define the function
\begin{align}
    \delta:\mathbb{R}^m_{>0}\rightarrow \mathbb{R},\ \ \delta(\boldsymbol{\alpha})=\sup \left\{ \frac{h(\mu)}{\lambda(\mu)}: \mu \in M(\tilde f),\lambda(\mu)<\infty,\ \int a_i d \mu = \alpha_i\ 1\leq i\leq m 
   \right\}.\nonumber
\end{align}

Using the following lemma, we can improve the conditional variational principle with $\epsilon$ in Proposition \ref{conditional variational principle with epsilon} to the conditional variational principle.

\begin{lemma}\label{continuity of delta}
The function $\delta$ is continuous.    
\end{lemma}
\begin{pf}
Let $\boldsymbol{\alpha}=(\alpha_1,\cdots,\alpha_m)\in \text{Int}(\text{Im}\Psi)$. Then, we can find measures $\mu_{\omega}\in\lbrace \mu \in \mathcal{M}({\tilde f}): \lambda(\mu)<\infty \rbrace \ (\omega=(\omega_1,\cdots,\omega_m)\in\lbrace1,2\rbrace^m)$ such that for $i=1,2,\cdots,m$,
\begin{align}
    \int a_i d\mu_\omega<\alpha_i\ \text{if}\ \omega_i=1,\ \text{and}\ 
    \int a_i d\mu_\omega>\alpha_i\ \text{if}\ \omega_i=2.\nonumber
\end{align}
Thus, we have $\boldsymbol{\alpha}\in D:=\text{Int}\left(\text{Conv}\left\{\Psi(\mu_{\omega}):\omega\in\lbrace1,2\rbrace^m\right\}\right)$.  

We shall show that the function $\delta$ is continuous on $D$. To do this, we consider a sequence $(\boldsymbol{\beta}_n)_{n\in\mathbb{N}}\subset D$ such that $\lim_{n\to\infty}\boldsymbol{\beta}_n=\boldsymbol{\alpha}$. First, we will show that 
$\delta(\boldsymbol{\alpha})\leq\limsup_{n\to\infty}\delta(\boldsymbol{\beta}_n)$.

Let $\epsilon>0$. By definition of $\delta$ and Lemma \ref{imformation of U} there exists $ \mu \in {M}({\tilde f})$ with $\lambda(\mu)<\infty$ such that $h(\mu)/\lambda(\mu)>\delta(\boldsymbol{\alpha})-\epsilon$ and $\Psi(\mu)=\boldsymbol{\alpha}$. Since $(\boldsymbol{\beta}_n)_{n\in\mathbb{N}}\subset D$ and $\lim_{n\to\infty}\boldsymbol{\beta}_n=\boldsymbol{\alpha}\in D$, there exists a sequence  
of probability vectors $(p_{n,1},\cdots,p_{n,2^m})_{n\in \mathbb{N}}$ with $\lim_{n\to\infty} (p_{n,1},\cdots,p_{n,2^m})= (0,\dots ,0,1)$ and a sequence of vectors  $(\omega^{n,1},\cdots,\omega^{n,2^m-1})_{n\in \mathbb{N}}$ with values in $\lbrace1,2\rbrace^m$ 
 such that
\begin{align}   \boldsymbol{\beta}_n=\sum_{i=1}^{2^m-1}p_{n,i}\Psi(\mu_{\omega^{n,i}})+p_{n,2^m}\boldsymbol{\alpha}.\nonumber
\end{align}

Define $\nu_n\in\lbrace \mu \in {M}({\tilde f}): \lambda(\mu)<\infty \rbrace$ by $\nu_n:=\sum_{i=1}^{2^m-1}p_{n,i}\mu_{\omega^{n,i}}+p_{n,2^m}\mu$. Since  
$\Psi(\nu_n)=\boldsymbol{\beta}_n$, 

\begin{align}
    \delta(\boldsymbol{\alpha})-\epsilon<\frac{h(\mu)}{\lambda(\mu)}
    =\lim_{n\to\infty}\frac{\sum_{i=1}^{2^m-1}p_{n,i}h(\mu_{\omega^{n,i}})+p_{n,2^m}h(\mu)}{\sum_{i=1}^{2^m-1}p_{n,i}\lambda(\mu_{\omega^{n,i}})+p_{n,2^m}\lambda(\mu)}\nonumber=\lim_{n\to\infty}\frac{h(\nu_n)}{\lambda(\nu_n)}\leq\liminf_{n\to\infty}\delta(\boldsymbol{\beta}_n).\nonumber
\end{align}
Letting $\epsilon\rightarrow0$, we get $\delta(\boldsymbol{\alpha})\leq\liminf_{n\to\infty}\delta(\boldsymbol{\beta}_n)$.

Next, we shall show that $\limsup_{n\to\infty}\delta(\boldsymbol{\beta}_n)\leq\delta(\boldsymbol{\alpha})$. For a contradiction, suppose that there exists $\eta$ such that $\delta(\boldsymbol{\alpha})<\eta<\limsup_{n\to\infty}\delta(\boldsymbol{\beta}_n)$.
Then, we have $\eta<\sup\lbrace\delta(\boldsymbol{\beta}_k):k\geq n \rbrace$ for all $n\in\mathbb{N}$. Therefore, there exists $k_1\geq1$ and $\mu_{k_1}\in\lbrace \mu \in \mathcal{M}_{\tilde f}: \lambda(\mu)<\infty \rbrace$ such that $\eta<h(\mu_{k_1})/\lambda(\mu_{k_1})$, and $\Psi(\mu_{k_1})=\boldsymbol{\beta}_{k_1}$. Then, there exists $k_2\geq k_1+1$ and $\mu_{k_2}\in\lbrace \mu \in \mathcal{M}_{\tilde f}: \lambda(\mu)<\infty \rbrace$ such that $\eta<h(\mu_{k_2})/\lambda(\mu_{k_2})$, and$\Psi(\mu_{k_2})=\boldsymbol{\beta}_{k_2}$. Repeating this procedure, we can find a sequence $(k_n)_{n\in\mathbb{N}}$ with $k_1<k_2<\cdots<k_n<\cdots$
and $(\mu_{k_n})_{n\in\mathbb{N}}\subset \lbrace \mu \in \mathcal{M}_{\tilde f}: \lambda(\mu)<\infty \rbrace$ such that
$\eta<{h(\mu_{k_n})}/{\lambda(\mu_{k_n})},$ and
$\Psi(\mu_{k_n})=\boldsymbol{\beta}_{k_n}.\nonumber$
 Since $(\boldsymbol{\beta}_n)_{n\in\mathbb{N}}\subset D$ and $\lim_{n\to\infty}\boldsymbol{\beta}_n=\boldsymbol{\alpha}\in D$,  there exists a sequence of probability vectors $(p_{k_n,1},\cdots,p_{k_n,2^m})_{n\in \mathbb{N}}$ with $\lim_{n\to\infty} (p_{k_n,1},\cdots,p_{k_n,2^m})= (0,\dots ,0,1)$ and a sequence of vectors  $(\omega^{n,1},\cdots,\omega^{n,2^m-1})_{n\in \mathbb{N}}$ with values in $\lbrace1,2\rbrace^m$ such that
 \begin{align*}
     \boldsymbol{\alpha}=\sum_{i=1}^{2^m-1}p_{k_n,i}\Psi(\mu_{\omega^{n,i}})+p_{k_n,2^m}\boldsymbol{\beta}_{k_n}.
 \end{align*}

Define  $\nu_{k_n}\in\lbrace \mu \in {M}({\tilde f}): \lambda(\mu)<\infty \rbrace$ by $\nu_{k_n}:=\sum_{i=1}^{2^m-1}p_{k_n,i}\mu_{\omega^{n,i}}+p_{k_n,2^m}\mu_{k_n}$. Since $\Psi(\nu_{k_n})=\boldsymbol{\alpha}$,
\begin{align}
    \delta(\boldsymbol{\alpha})<\eta\leq \liminf_{n\to\infty}\frac{h(\mu_{k_n})}{\lambda(\mu_{k_n})}
    =\liminf_{n\to\infty}\frac{\sum_{i=1}^{2^m-1}p_{k_n,i}h(\mu_{\omega^{n,i}})+p_{k_n,2^m}h(\mu_{k_n})}{\sum_{i=1}^{2^m-1}p_{k_n,i}\lambda(\mu_{\omega^{n,i}})+p_{k_n,2^m}\lambda(\mu_{k_n})}\nonumber
    =\liminf_{n\to\infty}\frac{h(\nu_{k_n})}{\lambda(\nu_{k_n})}\leq\delta(\boldsymbol{\alpha}).\nonumber
\end{align}

This is a contradiction. Hence, we get $\limsup_{n\to\infty}\delta(\boldsymbol{\beta}_n)\leq\delta(\boldsymbol{\alpha})$ and $\delta(\boldsymbol{\alpha})=\lim_{n\to\infty}\delta(\boldsymbol{\beta}_n)$. This means that $\delta$ is continuous on $\text{Int}(\text{Im}\Psi)$.
$\qed$  
\end{pf}

\begin{prop}\label{improve}
    For all $\boldsymbol{\alpha}\in\mathbb{R}^m_{>0}$ we have
    \begin{align}\label{eq;variatnal}
        b(\boldsymbol{\alpha})=\delta(\boldsymbol{\alpha})=\sup\left\{\frac{h(\mu)}{\lambda(\mu)}:\mu\in M(\tilde f),\ \lambda(\mu)<\infty,\ \Psi(\mu)=\boldsymbol{\alpha}\right\}.
    \end{align}
\end{prop}

\begin{pf}
Let $\boldsymbol{\alpha}\in \text{Int}(\text{Im}\Psi)$. By Proposition \ref{conditional variational principle with epsilon}, we have $\delta(\boldsymbol{\alpha})\leq b(\boldsymbol{\alpha})$.
Moreover, there exists $\boldsymbol{\alpha}_n\in \text{Int}(\text{Im}\Psi)$ and $\mu_n\in\lbrace \mu \in {M}({\tilde f}): \lambda(\mu)<\infty \rbrace$ such that $\lim_{n\to\infty}\boldsymbol{\alpha}_n=\boldsymbol{\alpha}$, $\Psi(\mu_n)=\boldsymbol{\alpha}_{n}$, and $\lim_{n\to\infty}h(\mu_n)/\lambda(\mu_n)=b(\boldsymbol{\alpha})$. By continuity of $\delta$ (Lemma \ref{continuity of delta}), we obtain
\begin{align}
b(\boldsymbol{\alpha})=\lim_{n\to\infty}\frac{h(\mu_n)}{\lambda(\mu_n)}\leq\lim_{n\to\infty}\delta(\boldsymbol{\alpha}_n)=\delta(\boldsymbol{\alpha}).\nonumber
\end{align}
This implies that $b(\boldsymbol{\alpha})=\delta(\boldsymbol{\alpha})$.
$\qed$
\end{pf}

\section{Convergence of the topological pressure and Bowen's formula}\label{Convergence of the topological pressure and Bowen's formula}
 Define the $\mathbb{R}^m$-valued function
    $\Phi:\Lambda_c(G)\rightarrow \mathbb{R}^m$ by $\Phi(x):=(-a_1(x),-a_2(x),\cdots,-a_m(x)). \nonumber$
We will use the notations $\mathbb{R}_{\geq0}^m:=[0,\infty)^m$ and $\mathbb{R}_{\geq0}:=[0,\infty)$    
The purpose of this section is to determine the set on which the function $(\boldsymbol{\alpha},\boldsymbol{q},b)\in\mathbb{R}_{\geq 0}^m\times\mathbb{R}^m\times\mathbb{R}\mapsto p(\boldsymbol{\alpha},\boldsymbol{q},b)$ is finite (see Lemma \ref{convergence lemma}). Also, by the proof of Lemma \ref{convergence lemma}, we obtain the improved Bowen's formula in Theorem \ref{useful Bowen}. Moreover, we will show that a measure of maximal dimension is the equilibrium measure of the geometric potential $-\dim_H(\Lambda_c(G))\log|\tilde f'|$ and we obtain item $(2)$ of Theorem \ref{main} in Proposition \ref{prop item (2)}. 

We shall use the notation $a\ll b$ for two positive reals $a,b$ to indicate that there exists a constant $C>0$ such that $a/b\leq C$. If $a\ll b$ and $b\ll a$, we write $a\asymp b$.

We first prove the following two lemmas used to determine the set on which the function $(\boldsymbol{\alpha},\boldsymbol{q},b)\in\mathbb{R}_{\geq 0}^m\times\mathbb{R}^m\times\mathbb{R}\mapsto p(\boldsymbol{\alpha},\boldsymbol{q},b)\in\mathbb{R}\nonumber$ is finite. 
These two lemmas can be obtained by a technical calculation.

\begin{lemma}\label{almost ^2}
    There exists constant $W>1$ such that 
\begin{align}\nonumber
    -2t(\log  W + \log l)
    \leq \sup_{ [\gamma^lh]}(-t\log |(\gamma^l)'| \circ \pi)
    \leq -2t(-\log  W + \log l)\nonumber
\end{align}
for all $l\geq 1$, $t>0$, $\gamma\in \Gamma_0,$ and $h\in H_0$.
\end{lemma}

\begin{pf}
    The proof  follows from the calculations in the proof of \cite[Lemma 2.8]{bowen1979markov}. $\qed$
\end{pf}

\begin{lemma}\label{asymptotic behavior of log}
We have 
\begin{align}\nonumber
\lim_{x \to p_i}\frac{a_i(x)}{\log|\tilde f'(x)|}= \infty\ \ (i=1,2,\cdots,m).
\end{align}

\end{lemma}

\begin{pf}
This is the immediately consequence of Lemma \ref{almost ^2}. $\qed$
\end{pf}

Let $\mathcal{F}:=\{(\boldsymbol{q},b)\in\mathbb{R}^m\times\mathbb{R}_{\geq 0}:P(\langle\boldsymbol{q},\Phi\rangle-b\log|\tilde f'|)<\infty\}$ and let $\partial\mathbb{R}_{\geq 0}^m:=\lbrace\boldsymbol{q}\in\mathbb{R}^m_{\geq 0}:\exists i \in \lbrace1,2,\cdots,m\rbrace\ \text{s.t.\ }q_i=0 \rbrace$. 
The following lemma determines the set $\mathcal{F}$.

\begin{lemma}\label{convergence lemma}
    We have $\mathcal{F}=(\mathbb{R}_{>0}^m\times\mathbb{R}_{\geq 0})\cup(\partial\mathbb{R}_{\geq 0}^m\times\{b\in\mathbb{R}_{\geq 0}:b>1/2\})$.
\end{lemma}
\begin{pf}
    Note that for all $\boldsymbol{q}\in\mathbb{R}_{>0}^m$ and $x\in\Lambda_c(G)$ we have $\sum_{\omega \in \mathcal{A}}\exp \left(\sup_{\tau \in [\omega]}  \langle\boldsymbol{q},\Phi(\pi(\tau))\rangle \right)\nonumber
=2\#H_0\sum_{i=1}^m\sum_{l=1}^{\infty}\exp(-q_i(l-1))<\infty\nonumber$ and $\log|\tilde f'(x)|>0$. By Lemma \ref{original convergence lemma}, we obtain 
    $P(\langle \boldsymbol{q}, \Phi \rangle -b \log|\tilde f'|)\leq P(\langle\boldsymbol{q},\Phi\rangle) < \infty$
 for all $\boldsymbol{q}\in\mathbb{R}^m_{>0}$ and $b\geq0$. On the other hand, by Lemma \ref{asymptotic behavior of log}, for each $\boldsymbol{q}\notin\mathbb{R}_{\geq 0}^m$ with $q_i<0$ $(i=1,2,\cdots,m)$ and
 $M:={b}/{-q_i}$ there exists $l_M >0$ such that 
$-q_i a_i(\pi(\gamma_i^{l}hh\cdots))\geq -q_iM\log|\tilde f'(\pi(\gamma_i^{l}hh\cdots))|$
  for all $l\geq l_M$, $h\in H_0$. Hence, we obtain 
\begin{align}
&\sum_{\omega \in \mathcal{A}}\exp \left(\sup_{\tau \in [\omega]}  (\langle\boldsymbol{q},\Phi(\pi(\tau))\rangle-b\log|\tilde f'(\pi\circ\tau)|) \right)\nonumber 
 \geq  \sum_{l=l_M}^\infty\exp{(-q_ia_i(\omega_{l})-b \log |\tilde f'(\omega_{l})|})\geq\sum_{l=l_M}^\infty1=\infty\nonumber
\end{align}
for $b\geq0$. 
Therefore, by Lemma \ref{original convergence lemma}, we obtain 
$P(\langle \boldsymbol{q}, \Phi \rangle -b \log|\tilde f'|)=\infty$   
 for all $\boldsymbol{q}\notin\mathbb{R}_{\geq 0}^m$ and $b\geq0$.
 Moreover, by Lemma \ref{almost ^2}, there exists $K>1$ such that
 \begin{align}\label{proof of Bowen's formula}
     \frac{1}{K}\sum_{l=0}^\infty l^{-2b}\leq\sum_{\omega\in\mathcal{A}} \exp\left(\sup_{\tau\in[\omega]}(-b\log|\tilde f'\circ\pi|)\right)\leq K\sum_{l=0}^\infty l^{-2b}
 \end{align}
 for all $b\geq0$. This implies that
     $P(\langle\boldsymbol{q},\Phi\rangle-b\log|\tilde f'|)<\infty
 \ \text{if\ and\ only\ if}\  b>{1}/{2}$  for $\boldsymbol{q}\in\partial\mathbb{R}_{\geq 0}^m$.
 $\qed$
\end{pf}
By Lemma \ref{convergence lemma} and the  proof of \cite[Proposition 2.1.9]{mauldin2003graph}, we have  $\nonumber
    p(\boldsymbol{\alpha},\boldsymbol{q},b(\boldsymbol{\alpha}))=\infty$
for all $\boldsymbol{q} \in  \lbrace (x_1,x_2,\cdots,x_m)\in \mathbb{R}^m:\exists i \in\lbrace1,2,\cdots,m\rbrace\ \text{s.t.\ } x_i<0\rbrace$, and $\boldsymbol{\alpha}\in \mathbb{R}^m_{\geq 0}$.

Note that by Lemma \ref{original convergence lemma}, the inequality $(\ref{proof of Bowen's formula})$ implies that $\lim_{b\searrow 1/2}P(-b\log|\tilde f'|)=\infty$. Therefore, by Theorem \ref{Bowen}, we obtain the following useful formula.

\begin{thm}\label{useful Bowen}(Bowen's formula)
    We have $P(-\dim_H(\Lambda_c(G))\log|\tilde f'|)=0$.  
\end{thm}
By Theorem \ref{useful Bowen}, we obtain item (2) of Theorem \ref{main}.
\begin{prop}\label{prop item (2)}
    Let $\mu$ be the equilibrium state of the geometric potential $-\dim_H{\Lambda_{c}(G)}\log |\tilde f'|$. The measure $\mu$ is the maximal measure that is 
    \begin{align}\label{maximal measure}
        \dim_H(\mu)=\dim_H(\Lambda_c(G))=\dim_H(J(\infty,\infty,\cdots,\infty)).
    \end{align}
\end{prop}
\begin{pf}
   Let $\mu$ is the equilibrium measure of the geometric potential $-\dim_H{\Lambda_{c}(G)}\log |\tilde f'|$. Then, using Theorem \ref{useful Bowen}, we have 
$\dim_H(\Lambda_c(G))=\dim_H(\mu)$. 
Next, we shall show the second equality of (\ref{maximal measure}). Note that $\mu$ is a Gibbs measure. By Lemma \ref{almost ^2} and Theorem \ref{useful Bowen}, for all $1\leq i\leq m$ we have
\begin{align*}
    \int a_i d\mu&=\sum_{\gamma \in \{\gamma_i,\gamma_{i}^{-1}\},\ h\in H_0}\sum_{l=1}^{\infty}(l-1)\mu([\gamma^lh])\asymp\sum_{\gamma \in \{\gamma_i,\gamma_i^{-1}\},\ h\in H_0}\sum_{l=1}^{\infty}(l-1)\exp\left(\sup_{\tau\in[\gamma^lh]}(-b\log|\tilde f'\circ\pi|)\right)\\&\asymp\sum_{l=1}^\infty l^{-2\dim_H(\Lambda_c(G))+1}=\infty.
\end{align*}
Thus, by Proposition \ref{conditional variational principle with epsilon}, we get
\begin{align*}
    \dim_H(\Lambda_c(G))=\dim_H(\mu)=\frac{h(\mu)}{\lambda(\mu)}\leq \dim_H(J(\infty,\infty,\cdots,\infty))\leq\dim_H(\Lambda_c(G)).\ \qed
\end{align*}
\end{pf}

\section{Relationship between the topological pressure and dimension spectrum}\label{Relationship between the topological pressure and dimension spectrum}
In this section, we relate the topological pressure and the dimension spectrum. To do this, for all $\boldsymbol{\alpha}\in \text{Int}(\text{Im}\Psi)=\mathbb{R}_{>0}^m$ (see Lemma \ref{imformation of U}) we consider the graph of the function $\boldsymbol{q}\in\mathbb{R}^m_{\geq 0}\mapsto p(\boldsymbol{\alpha},\boldsymbol{q},b(\boldsymbol{\alpha}))\in\mathbb{R}$ and show in Proposition \ref{prop more important proposition}
  that there exists $\boldsymbol{q}(\boldsymbol{\alpha})\in\mathbb{R}_{>0}^m$ such that the function  $\boldsymbol{q}\in\mathbb{R}_{\geq 0}^m\mapsto p(\boldsymbol{\alpha},\boldsymbol{q},b(\boldsymbol{\alpha}))\in\mathbb{R}$ takes the minimum value $0$ at $\boldsymbol{q}(\boldsymbol{\alpha})$.  

Recall that $p(\boldsymbol{\alpha},\boldsymbol{q},b(\boldsymbol{\alpha}))=P(\langle\boldsymbol{q},(\Phi+\boldsymbol{\alpha})\rangle-b(\boldsymbol{\alpha})\log|\tilde f'|)<\infty$ for all $\boldsymbol{q},\boldsymbol{\alpha}\in\mathbb{R}_{>0}^m$ by Lemma \ref{convergence lemma}.
\begin{lemma}\label{positivity}
We have 
$\nonumber
     p(\boldsymbol{\alpha},\boldsymbol{q},b(\boldsymbol{\alpha}))\geq0$
for all $\boldsymbol{q}\in \mathbb{R}^m\  \text{and}\  \boldsymbol{\alpha}\in  \mathbb{R}^m_{>0}$.     
\end{lemma}

\begin{pf}
 Let $\boldsymbol{\alpha}=(\alpha_1,\alpha_2,\cdots,\alpha_m)\in \mathbb{R}^m_{>0}$. By Proposition \ref{improve} there exists   
a sequence of $\tilde f$-invariant measures $(\mu_n)_n$ such that the following conditions hold:
(1) $\int a_i d\mu_n = \alpha_i$ for $i\in\lbrace1,2,\cdots,m\rbrace$.
(2) $h(\mu_n)<\infty$ and $\lambda(\mu_n)<\infty$.  
(3) $\lim_{n \to \infty} {h(\mu_n)}/{\lambda(\mu_n)} = b(\boldsymbol{\alpha})$.
If we choose $0<s_1<s_2< b(\boldsymbol{\alpha})$ then for all $\boldsymbol{q}_0=(q_{0,1},q_{0,2},\cdots,q_{0,m})\in \mathbb{R}_{>0}^m$ we have $K:=P(\langle \boldsymbol{q_0},\Phi\rangle- s_1 \log |\tilde f'|) <\infty$ by Lemma \ref{convergence lemma}. Fix $\boldsymbol{q}_0\in\mathbb{R}_{>0}^m$. By the variational principle for the topological pressure, we have 
for all $n \in \mathbb{N}$, 
\begin{align}\label{bounded}
- s_1 \lambda(\mu_n)+h(\mu_n)\leq K +\langle \boldsymbol{q}_0, \boldsymbol{\alpha}\rangle.
\end{align}
On the other hand, since we have 
$s_2 < {h(\mu_n)}/{\lambda(\mu_n)} \leq b(\boldsymbol{\alpha})$ for $n$ sufficiently large, we obtain $s_2 \lambda(\mu_n) \leq h(\mu_n)$. Thus, we obtain $h(\mu_n)-s_1 \lambda(\mu_n) \geq (s_2 - s_1)\lambda(\mu_n)$ for $n$ sufficiently large. Substituting this into inequality ($\ref{bounded}$), we obtain 
$\nonumber
\lambda(\mu_n) \leq ({K+\langle\boldsymbol{q_0},\boldsymbol{\alpha}\rangle})/({s_2-s_1})$
  for $n$ sufficiently large. This means that $\sup_{n\in\mathbb{N}}\lambda(\mu_n)<\infty$. Furthermore, by the variational principle for the topological pressure, we have 
\begin{align}\nonumber
&p(\boldsymbol{\alpha},\boldsymbol{q},b(\boldsymbol{\alpha})) \nonumber
\geq  h(\mu_n) - b(\boldsymbol{\alpha}) \lambda(\mu_n)= \lambda(\mu_n) \left (\frac{h(\mu_n)}{ \lambda(\mu_n)} -b(\boldsymbol{\alpha}) \right )\nonumber
\end{align}
for all $\boldsymbol{q} \in \mathbb{R}^m$ and $\boldsymbol{\alpha}\in \mathbb{R}^m_{>0}$. Since  $\lambda(\mu_n)$ is bounded above,
 $\nonumber
 \lim_{n\to \infty}\lambda(\mu_n) \left ({h(\mu_n)}/{ \lambda(\mu_n)} -b(\boldsymbol{\alpha}) \right ) =0.$
 This means that 
 $p(\boldsymbol{\alpha},\boldsymbol{q},b(\boldsymbol{\alpha}))\geq0$
for all $\boldsymbol{q}\in \mathbb{R}^m\  \text{and}\  \boldsymbol{\alpha}\in \mathbb{R}^m_{>0}$.
$\qed$
\end{pf}

For $n\in\mathbb{N}$ we will use the notations $M_n:=\bigcup_{i=1}^n\lbrace\gamma^ih:\gamma\in\Gamma_0,h\in H_0\rbrace\cup H_0$ and $K_n:=M_n^{\mathbb{N}}\cap\Sigma_A$. For a weakly H\"older potential $\phi:\Lambda_c(G)\rightarrow\mathbb{R}$ and $n\in\mathbb{N}$ we define
\begin{align*}
    P_{\pi(K_n)}(\phi)
    :=\sup \left\{ h(\mu) + \int \phi d\mu : \mu \in {M}({\tilde f,\pi(K_n)}) \right\},
\end{align*}
where ${M}({\tilde f,\pi(K_n)})$ denotes the set of $\tilde f$-invariant Borel probability measure supported on $\pi(K_n)$. Note that $\pi:K_n\rightarrow \pi(K_n)$ is homeomorphism. Hence, by the variational principle \cite[Theorem 3.4.1]{przytycki2010conformal}, for a weakly H\"older potential $\phi:\Lambda_c(G)\rightarrow\mathbb{R}$ and $n\in\mathbb{N}$ we obtain $P_{\pi(K_n)}(\phi)=P_{M_n}(\phi)$ (see Definition \ref{Definition of pressure} for the definition of $P_{M_n}(\phi)$). For simplicity 
we will use the notation $p_n(\boldsymbol{q},\boldsymbol{\alpha},b(\boldsymbol{\alpha}))=P_{\pi(K_n)}(\langle \boldsymbol{q}, (\Phi+\boldsymbol{\alpha}) \rangle -b(\boldsymbol{\alpha}) \log|\tilde f'|)$ for $\boldsymbol{\alpha},\boldsymbol{q}\in\mathbb{R}^m$ and $n\in\mathbb{N}$.

The following lemma will be used in the proof of Proposition \ref{important lemma}. 
Our case is in higher dimensions, but the proof of the following lemma is essentially the same as the proof of \cite[Lemma 3.2]{iommi2017transience}.

\begin{lemma}\label{positive approximation}
  If $\boldsymbol{\alpha}\in\mathbb{R}_{>0}^m$ and 
    $\inf \lbrace p(\boldsymbol{\alpha},\boldsymbol{q},b(\boldsymbol{\alpha})):\boldsymbol{q}\in\mathbb{R}^m\rbrace>0$
    then there exists $N\in\mathbb{N}$ such that
   $(1)\,\,\, p_{N}(\boldsymbol{q},\boldsymbol{\alpha},b(\boldsymbol{\alpha}))>0\ \text{for\ all}\ \boldsymbol{q}\in\mathbb{R}^m\ \ \textrm{ and }
(2)\,\,\, \lim_{|\boldsymbol{q}|\to\infty}p_{N}(\boldsymbol{q},\boldsymbol{\alpha},b(\boldsymbol{\alpha}))=\infty$ \textrm{ hold}. 
\end{lemma}
\begin{pf}
    Let $\boldsymbol{\alpha}=(\alpha_1,\cdots,\alpha_m)\in \mathbb{R}^m_{>0}$. 
    We start with the second part. We take $C\in\mathbb{N}$ such that
    \begin{align*}
        C> 2(m+1)\cdot{\max\lbrace |\alpha_k|:1\leq k\leq m\rbrace}\ \text{and\ put\ }  x_{i}^{\omega}=\left\{
 \begin{array}{ll}
   \gamma_i^Ch_1   & \text{if}\ \omega=1 \\
   h_1   & \text{if} \ \omega=2\nonumber
 \end{array}
 \right. 
    \end{align*}
    {for} $\omega\in\lbrace1,2\rbrace$ and $i\in\{1,2,\cdots,m\}$.
 For $\tau\in\lbrace1,2\rbrace^m$ we define    
  \begin{align*}
      &x_{\tau}:=\pi(x_{1}^{\tau_1}h_1x_{2}^{\tau_2}h_1\cdots x_{m}^{\tau_m}h_1x_{1}^{\tau_1}h_1x_{2}^{\tau_2}h_1\cdots x_{m}^{\tau_m}h_1\cdots).
  \end{align*}   
For $\tau\in\lbrace1,2\rbrace^m$ we consider $\tilde f$-invariant measures $\mu_{\tau}:=1/2m\sum_{i=0}^{2m-1}\delta_{\tilde f^i(x_{\tau})}$ supported on $\pi(K_C)$.    
For $\tau\in\lbrace1,2\rbrace^m$ and $i\in\{1,2,\cdots,m\}$ we have
\begin{align}
    \int a_i d\mu_\tau>\alpha_i\ \text{if}\ \tau_i=1,\ \text{and}\ 
    \int a_i d\mu_\tau=0<\alpha_i\ \text{if}\ \tau_i=2\nonumber.
\end{align}
Let $(\boldsymbol{q}_n)_{n\in\mathbb{N}}=(q_{n,1},q_{n,2},\cdots,q_{n,m})_{n\in\mathbb{N}}\subset\mathbb{R}^m_{\geq 0}$ with $\lim_{n\to\infty}|\boldsymbol{q}_n|=\infty$.
First assume that there exists $i\in \lbrace1,2,\cdots,m\rbrace$ such that $\lim_{n\to\infty}q_{n,i}= \infty$. Take $\tau\in\{1,2\}^m$ with $\tau_i=2$. By variational principle for the topological pressure, we have
\begin{align}
    \lim_{n \to \infty}p_{C}(\boldsymbol{\alpha},\boldsymbol{q}_n,b(\boldsymbol{\alpha}))\nonumber\nonumber
    \geq \lim_{n\to\infty}q_{n,i}\left(\int -a_id\mu_{\tau} +\alpha_i\right)-b(\boldsymbol{\alpha})\lambda(\mu_{\tau})=\infty.\nonumber
\end{align} 
The remaining case is similar and therefore omitted.

Next, we shall show item $(1)$. Suppose for a contradiction that $(1)$ does not hold. Then, for each $n\in\mathbb{N}$ with $n\geq C$ there exists $\boldsymbol{q}_n\in\mathbb{R}^m$ such that $p_{n}(\boldsymbol{\alpha},\boldsymbol{q}_n,b(\boldsymbol{\alpha}))\leq0$. By $(2)$, we have $\sup_{n\geq C}|\boldsymbol{q}_n|<\infty$. Hence, by passing to subsequence, we may assume that there exists $\boldsymbol{q}_*\in\mathbb{R}^m$ such that $\lim_{n\to\infty}\boldsymbol{q}_n=\boldsymbol{q}_*$. By Proposition \ref{approximation property}, we obtain
\begin{align*}
    p(\boldsymbol{\alpha},\boldsymbol{q}_*,b(\boldsymbol{\alpha}))=\lim_{n\to\infty}p_{n}(\boldsymbol{\alpha},\boldsymbol{q}_*,b(\boldsymbol{\alpha}))
    =\lim_{n\to\infty}\lim_{k\to\infty}p_{n}(\boldsymbol{\alpha},\boldsymbol{q}_k,b(\boldsymbol{\alpha}))\leq\limsup_{n\to\infty}p_{n}(\boldsymbol{\alpha},\boldsymbol{q}_n,b(\boldsymbol{\alpha}))\leq0.
\end{align*}
This is a contradiction.
$\qed$

\end{pf}

   By Lemma \ref{convergence lemma} and \cite[Proposition 2.6.12]{mauldin2003graph}, we obtain the following formula for   the derivative of the  topological pressure which is called Ruelle's formula: For all $\boldsymbol{\alpha},\boldsymbol{q}\in\mathbb{R}^m_{>0}$ and $i\in\{1,2,\cdots,m\}$ 
   \begin{align}\label{ruelle's formula first}
   \frac{\partial}{\partial q_i}p(\boldsymbol{\alpha},\boldsymbol{q},b(\boldsymbol{\alpha}))=\int(-a_i+\alpha_i)d\mu_{\boldsymbol{q}}
   \end{align}
   where $\mu_{\boldsymbol{q}}$ denote the equilibrium state of the potential $\langle\boldsymbol{q},\Phi+\boldsymbol{\alpha}\rangle-b(\boldsymbol{\alpha})\log|\tilde f'|)$.

   The following lemma state that there are two types of relation between the topological pressure and the dimension spectrum. One type is that the function $\boldsymbol{q}\in\mathbb{R}^m_{\geq 0}\mapsto p(\boldsymbol{\alpha},\boldsymbol{q},b(\boldsymbol{\alpha}))\in\mathbb{R}$ attains its minimum value on the boundary of $\mathbb{R}_{\geq 0}^m$. The other type is that the function $\boldsymbol{q}\in\mathbb{R}^m_{\geq 0}\mapsto p(\boldsymbol{\alpha},\boldsymbol{q},b(\boldsymbol{\alpha}))\in\mathbb{R}$ attains its minimum value interior of $\mathbb{R}_{\geq 0}^m$. As mentioned in the introduction, the latter type is the desired type. Furthermore, we will show later that the later type always occurs.

\begin{prop}\label{important lemma}
For all $\boldsymbol{\alpha}\in \mathbb{R}^m_{>0}$ there exists $\boldsymbol{q}(\boldsymbol{\alpha})\in\mathbb{R}_{\geq 0}^m$ such that $p(\boldsymbol{\alpha},\boldsymbol{q}(\boldsymbol{\alpha}),b(\boldsymbol{\alpha}))=0$.
\end{prop}

\begin{pf}
Let $\boldsymbol{\alpha}\in \mathbb{R}^m_{>0}$. By virtue of Lemma \ref{positivity}, we have 
    $p(\boldsymbol{\alpha},\boldsymbol{q},b(\boldsymbol{\alpha}))  \geq 0$
 for all $\boldsymbol{q} \in \mathbb{R}^m$. We first assume that 
   there exists $\boldsymbol{q}(\boldsymbol{\alpha})=(q_{\boldsymbol{\alpha},1},q_{\boldsymbol{\alpha},2},\cdots,q_{\boldsymbol{\alpha},m}) \in \mathbb{R}^m_{>0}$ such that
     ${\partial}/{\partial q_i}p(\boldsymbol{\alpha},\boldsymbol{q}(\boldsymbol{\alpha}),b(\boldsymbol{\alpha}))=0$ $(i=1,2,\cdots,m)$. We will show that 
     $p(\boldsymbol{\alpha},\boldsymbol{q}(\boldsymbol{\alpha}),b(\boldsymbol{\alpha}))=0.$
  Note that $p(\boldsymbol{\alpha},\boldsymbol{q}(\boldsymbol{\alpha}),b(\boldsymbol{\alpha}))$ is finite by Lemma \ref{convergence lemma}, and $-a_i\ (i=1,2,\cdots,m)$ and $\log |\tilde f'|$ are weakly H\"older. Thus, there exists an equilibrium measure $\mu_{\boldsymbol{q}(\boldsymbol{\alpha})}$ for the potential $\langle \boldsymbol{q}(\boldsymbol{\alpha}), (\Phi+\boldsymbol{\alpha}) \rangle -b(\boldsymbol{\alpha})\log|\tilde f'|$. By Ruelle's formula for the derivative of pressure (\ref{ruelle's formula first}), for $i\in \lbrace1,2,\cdots,m\rbrace$ we have
\begin{align}\nonumber
     \frac{\partial}{\partial q_i}p(\boldsymbol{\alpha},\boldsymbol{q}(\boldsymbol{\alpha}),b(\boldsymbol{\alpha}))=\int-a_id\mu_{\boldsymbol{q}(\boldsymbol{\alpha})}+\alpha_i=0 
 \end{align}
 and thus,  
 $\int a_i d\mu_{\boldsymbol{q}(\boldsymbol{\alpha})}=\alpha_i.$
 Hence, we obtain 
\begin{align}\nonumber
0\leq p(\boldsymbol{\alpha},\boldsymbol{q}(\boldsymbol{\alpha}),b(\boldsymbol{\alpha}))= -b(\boldsymbol{\alpha})\lambda(\mu_{\boldsymbol{q}(\boldsymbol{\alpha})})+h(\mu_{\boldsymbol{q}(\boldsymbol{\alpha})})
\leq -\frac{h(\mu_{\boldsymbol{q}(\boldsymbol{\alpha})})}{\lambda(\mu_{\boldsymbol{q}(\boldsymbol{\alpha})})}\lambda(\mu_{\boldsymbol{q}(\boldsymbol{\alpha})}) + h(\mu_{\boldsymbol{q}(\boldsymbol{\alpha})}) = 0
\end{align}
Thus, we conclude that 
\begin{align}\label{eq:attain measure}
p(\boldsymbol{\alpha},\boldsymbol{q}(\boldsymbol{\alpha}),b(\boldsymbol{\alpha}))=0\ \text{and}\  b(\boldsymbol{\alpha})=\frac{h(\mu_{\boldsymbol{q}(\boldsymbol{\alpha})})}{\lambda(\mu_{\boldsymbol{q}(\boldsymbol{\alpha})})}.
\end{align}
By Proposition \ref{improve}, the equation $b(\boldsymbol{\alpha})={h(\mu_{\boldsymbol{q}(\boldsymbol{\alpha})})}/{\lambda(\mu_{\boldsymbol{q}(\boldsymbol{\alpha})})}$ implies that $\mu_{\boldsymbol{q}(\boldsymbol{\alpha})}$ attains the spremum of the conditional variational principle $(\ref{eq;variatnal})$ for $\boldsymbol{\alpha}$.

Next, we assume that there is no point $\boldsymbol{q}_0=(q_{0,1},q_{0,2},\cdots,q_{0,m}) \in \mathbb{R}^m_{>0}$ such that
     ${\partial}/{\partial q_i}p(\boldsymbol{\alpha},\boldsymbol{q}(\boldsymbol{\alpha}),b(\boldsymbol{\alpha}))=0$
for $i\in\lbrace1,2,\cdots,m\rbrace$.
For a contradiction we assume that there exists $\boldsymbol{q}_0\in\partial\mathbb{R}_{\geq 0}^m$ such that 
$p(\boldsymbol{\alpha},\boldsymbol{q}_0,b(\boldsymbol{\alpha}))=\infty$.
Then, by Lemma \ref{convergence lemma}, we have 
    $p(\boldsymbol{\alpha},\boldsymbol{q},b(\boldsymbol{\alpha}))=\infty$
for all $\boldsymbol{q}\in\partial \mathbb{R}_{\geq 0}^m$.
Considering compact approximations of the topological pressure, we can show that the function $\boldsymbol{q}\in\mathbb{R}_{\geq0}^m\rightarrow p(\boldsymbol{\alpha},\boldsymbol{q},b(\boldsymbol{\alpha}))\in\mathbb{R}$ is lower semicontinuous and therefore, for all $\boldsymbol{q}\in \partial \mathbb{R}_{\geq 0}^m$ and a sequence $(\boldsymbol{q}_n)_{n\in\mathbb{N}}\subset\mathbb{R}_{>0}^m$ with $\lim_{n\to\infty}\boldsymbol{q}_{n}=\boldsymbol{q}$ we have
\begin{align*}
\lim_{\boldsymbol{q}_n \to \boldsymbol{q} } p(\boldsymbol{\alpha},\boldsymbol{q}_n,b(\boldsymbol{\alpha}))=\infty.
\end{align*}
By the same argument as in the proof of Lemma \ref{positive approximation}, for each sequence $(\boldsymbol{q}_n)_{n\in\mathbb{N}}\subset\mathbb{R}^m_{\geq 0}$ with $\lim_{n\to\infty}|\boldsymbol{q}_n|=\infty$ we have 
    $\lim_{n \to \infty}p(\boldsymbol{\alpha},\boldsymbol{q}_n,b(\boldsymbol{\alpha}))=\infty.$
Hence, there exists $\boldsymbol{q}_0=(q_{0,1},q_{0,2},\cdots,q_{0,m}) \in \mathbb{R}^m_{>0}$ such that the function $\boldsymbol{q}\in\mathbb{R}_{\geq0}^m\rightarrow p(\boldsymbol{\alpha},\boldsymbol{q},b(\boldsymbol{\alpha}))\in\mathbb{R}$  attains its minimum at $\boldsymbol{q}_0$ and thus,
$\nonumber
     {\partial}/{\partial q_i}p(\boldsymbol{\alpha},\boldsymbol{q}_0,b(\boldsymbol{q}(\boldsymbol{\alpha})))=0$
for $i\in\lbrace1,2,\cdots,m\rbrace$. This is a contradiction.

Also, for a contradiction, we assume that 
    $0<p(\boldsymbol{\alpha},\boldsymbol{q},b(\boldsymbol{\alpha}))<\infty$
for all $\boldsymbol{q}\in\partial\mathbb{R}_{\geq 0}^m$. Then, we have 
\begin{align*}
       0<p(\boldsymbol{\alpha},\boldsymbol{q},b(\boldsymbol{\alpha}))<\infty
       \end{align*}
for all $\boldsymbol{q}\in\mathbb{R}_{\geq 0}^m$.
By Lemma \ref{positive approximation}, there exists a compact set $N\in\mathbb{N}$ such that the function  $\boldsymbol{q}\in\mathbb{R}^m\mapsto p_N(\boldsymbol{\alpha},\boldsymbol{q},b(\boldsymbol{\alpha}))$ satisfies $(1)$ and $(2)$ of Lemma \ref{positive approximation}. Hence, there exists $\boldsymbol{q}_0=(q_{0,1},q_{0,2}\cdots,q_{0,m})\in\mathbb{R}^m$ such that the function  $\boldsymbol{q}\in\mathbb{R}^m\mapsto p_N(\boldsymbol{\alpha},\boldsymbol{q},b(\boldsymbol{\alpha}))$ attains its minimum at $\boldsymbol{q}_0\in\mathbb{R}^m$ and thus, 
\begin{align}
    \left. \frac{\partial}{\partial q_i}p_N(\boldsymbol{\alpha},\boldsymbol{q},b(\boldsymbol{\alpha}))\right|_{q_i=q_{0,i}}= 0\nonumber
\end{align}
for all $i\in\lbrace1,\cdots,m\rbrace$. Denote by $\mu_N$ the equilibrium measure for $\langle\boldsymbol{q}_0,\Phi+\boldsymbol{\alpha}\rangle-b(\boldsymbol{\alpha})\log|\tilde f '|$. We obtain 
 $\int a_i d \mu_{q_N} = \alpha_i$ for all $i\in\lbrace1,\cdots,m\rbrace$. Thus, we obtain 
\begin{align}\nonumber
0<p_N(\boldsymbol{\alpha},\boldsymbol{q}_0,b(\boldsymbol{\alpha}))=h(\mu_N)-b(\boldsymbol{\alpha})\lambda(\mu_N).\nonumber
\end{align}
and ${h(\mu_N)}/{\lambda(\mu_N)} > b(\boldsymbol{\alpha})$. This contradicts Proposition \ref{improve}. Hence, there exists $\boldsymbol{q}(\boldsymbol{\alpha})\in\partial\mathbb{R}_{\geq 0}^m$ such that $p(\boldsymbol{\alpha},\boldsymbol{q}(\boldsymbol{\alpha}),b(\boldsymbol{\alpha}))=0$.
$\qed$
\end{pf}
\begin{prop}\label{prop more important proposition}
    For all $\boldsymbol{\alpha}\in\mathbb{R}_{>0}^m$ there exists $\boldsymbol{q}(\boldsymbol{\alpha})\in\mathbb{R}_{>0}^m$ such that 
\begin{align}\nonumber
     p(\boldsymbol{\alpha},\boldsymbol{q}(\boldsymbol{\alpha}),b(\boldsymbol{\alpha}))=0\ \ \ 
        \text{and}\ \ \ 
\frac{\partial}{\partial q_i}p(\boldsymbol{\alpha},\boldsymbol{q}(\boldsymbol{\alpha}),b(\boldsymbol{\alpha}))=0\ \ (i=1,2,\cdots,m).
\end{align}
Moreover, we obtain 
\begin{align*}
    b(\boldsymbol{\alpha})=\frac{h(\mu_{\boldsymbol{q}(\boldsymbol{\alpha})})}{\lambda(\mu_{\boldsymbol{q}(\boldsymbol{\alpha})})}=\max\left\{\frac{h(\mu)}{\lambda(\mu)}:\mu\in M(\tilde f),\ \lambda(\mu)<\infty,\ \Psi(\mu)=\boldsymbol{\alpha}\right\}
    \end{align*}
where $\mu_{\boldsymbol{q}(\boldsymbol{\alpha})}$ denotes the equilibrium state of the potential $\langle\boldsymbol{q}(\boldsymbol{\alpha}),\Phi+\boldsymbol{\alpha}\rangle-b(\boldsymbol{\alpha})\log|\tilde f'|$
\end{prop}

\begin{pf}
 By Proposition \ref{important lemma}, it is suffices prove that the minimum of the function $\boldsymbol{q}\in\mathbb{R}^m_{\geq 0}\mapsto p(\boldsymbol{\alpha},\boldsymbol{q},b(\boldsymbol{\alpha}))\in\mathbb{R}$ is not attained at $\partial{R_{+}^m}$. For a contradiction, we assume that there exists $\boldsymbol{q}(\boldsymbol{\alpha})=(q_{\boldsymbol{\alpha},1},q_{\boldsymbol{\alpha},2},\cdots,q_{\boldsymbol{\alpha},m})\in\partial\mathbb{R}_{\geq 0}^m$ and $i\in\{1,2,\cdots,m\}$ such that $p(\boldsymbol{\alpha},\boldsymbol{q}(\boldsymbol{\alpha}),b(\boldsymbol{\alpha}))=0$ and $q_{\boldsymbol{\alpha},i}=0$. We consider a sequence $(\boldsymbol{q}_n)_{n\in\mathbb{N}}=((q_{n,1},\cdots,q_{n,m}))_{n\in\mathbb{N}}$ such that 
 \begin{align*}
 q_{n,j}=
 \left\{
 \begin{array}{ll}
   q_{\boldsymbol{\alpha},j}   & \text{if}\ j\neq i \\
   {1}/{n}   & \text{if} \ j=i\nonumber
 \end{array}
 \right.  
 \end{align*}
 for all $n\in\mathbb{N}$ and $j\in\{1,2,\cdots,m\}$.
    By convexity of the pressure functional, the following one-sided derivative exists: 
    \begin{align}
        \left. \frac{\partial}{\partial q_i}p(\boldsymbol{\alpha},q_{\boldsymbol{\alpha},1},\cdots,q_i,\cdots,q_{\boldsymbol{\alpha},m},b(\boldsymbol{\alpha}))\right|_{q_i=0_+}\nonumber
        =\lim_{n\to\infty}\left. \frac{\partial}{\partial q_i}p(\boldsymbol{\alpha},q_{\boldsymbol{\alpha},1},\cdots,q_i,\cdots,q_{\boldsymbol{\alpha},m},b(\boldsymbol{\alpha}))\right|_{q_i=\frac{1}{n}}.\nonumber
    \end{align}
    
    For all $n\in\mathbb{N}$ there exists a equilibrium measure $\nu_n$ for the potential $(\sum_{j\neq i,\ j\in\lbrace1,2,\cdots,m\rbrace}q_{\boldsymbol{\alpha},j}(-a_j+\alpha_j))+1/n(-a_i+\alpha_i)-b(\boldsymbol{\alpha})\log|\tilde f'|$. Note that $\nu_n$ is a Gibbs measure and  $0\leq\sup_{n\in\mathbb{N}}p(\boldsymbol{\alpha},\boldsymbol{q}_n,b(\boldsymbol{\alpha}))<\infty$. By Lemma \ref{almost ^2} and (\ref{ruelle's formula first}), we obtain 
\begin{align}
&\left. \frac{\partial}{\partial q_i}p(\boldsymbol{\alpha},q_{\boldsymbol{\alpha},1},\cdots,q_i,\cdots,q_{\boldsymbol{\alpha},m},b(\boldsymbol{\alpha}))\right|_{q_i=\frac{1}{n}}\nonumber
=\int(-a_i+\alpha_i)d\nu_{n}
\asymp\sum_{\gamma\in\lbrace\gamma_i,\gamma_i^{-1}\rbrace,\ h\in H_0}\sum_{l= 1}^{\infty}\int_{[\gamma^lh]} -a_i d\nu_{n}+\alpha_i \nonumber \\
&=\sum_{\gamma\in\lbrace\gamma_i,\gamma_i^{-1}\rbrace,\ h\in H_0}\sum_{l= 1}^{\infty} -l \nu_{n}(\pi([\gamma^{l+1}h]))\nonumber+\alpha_i\\ \nonumber
&\asymp  \sum_{\gamma\in\lbrace\gamma_i,\gamma_i^{-1}\rbrace,\ h\in H_0}\sum_{l= 1}^{\infty}  -l\cdot \exp \left(-p(\boldsymbol{\alpha},\boldsymbol{q}_n,b(\boldsymbol{\alpha}))-\frac{l}{n} +\frac{\alpha_i}{n}\right)\cdot \sup_{\pi([\gamma^{l+1}h])} |\tilde f'|^{-b(\boldsymbol{\alpha})}+\alpha_i  \\   \nonumber
&\asymp  \sum_{\gamma\in\lbrace\gamma_i,\gamma_i^{-1}\rbrace,\ h\in H_0}\sum_{l= 1}^{\infty}  -l\cdot \exp\left(-\frac{l}{n} \right)\cdot l^{-2b(\boldsymbol{\alpha})}  +\alpha_i\nonumber
\asymp  \sum_{l= 1}^{\infty} -\exp\left(-\frac{l}{n}\right) \cdot l^{-2b(\boldsymbol{\alpha})+1}+\alpha_i
\end{align}
 On the other hand, since $0\leq b(\boldsymbol{\alpha})\leq1$, we have
\begin{align}
\lim_{n \to \infty} \sum_{l= 1}^{\infty} - \exp\left(-\frac{l}{n}\right) \cdot l^{-2b(\boldsymbol{\alpha})+1}=
&\inf_{n\in\mathbb{N}}\sum_{l= 1}^{\infty}-\exp\left(-\frac{l}{n}\right) \cdot l^{-2b(\boldsymbol{\alpha})+1}=-\infty \nonumber.
\end{align}
Thus, we obtain 
    \begin{align}
        \left. \frac{\partial}{\partial q_i}p(\boldsymbol{\alpha},q_{\boldsymbol{\alpha},1},\cdots,q_i,\cdots,q_{\boldsymbol{\alpha},m},b(\boldsymbol{\alpha}))\right|_{q_i=0_+}\nonumber
        =\lim_{n\to\infty}\left. \frac{\partial}{\partial q_i}p(\boldsymbol{\alpha},q_{\boldsymbol{\alpha},1},\cdots,q_i,\cdots,q_{\boldsymbol{\alpha},m},b(\boldsymbol{\alpha}))\right|_{q_i=\frac{1}{n}}<0.\nonumber
    \end{align}

    This implies that there exists a neighborhood $\mathcal{V}\subset \mathbb{R}^m$ of $\boldsymbol{q}(\boldsymbol{\alpha})$ such that
$p(\boldsymbol{\alpha},\boldsymbol{q},b(\boldsymbol{\alpha}))<0$ for some $\boldsymbol{q}\in\mathcal{V}\cap\mathbb{R}_{>0}^m$. This is a contradiction. The proof is complete.
$\qed$
\end{pf}

\section{Regularity of the dimension spectrum}\label{Regularity of the dimension spectrum on a open "good" set}
In this section, we prove that the dimension spectrum is real-analytic on $\mathbb{R}_{>0}^m$ in Proposition \ref{prop regularity of the dimension spectrum}. Moreover, we prove item (4) of Theorem \ref{main} (see Proposition \ref{irregular set}). 

For a topological space $X$ and dynamical system $T:X\rightarrow X$ two continuous function $g_1,g_2:X\rightarrow \mathbb{R}$ are cohomologous if there exists a bounded continuous function $u:X\rightarrow\mathbb{R}$ and constant $c\in\mathbb{R}$ such that $g_1-g_2=u-u\circ T+c$. It is not difficult to verify that this relation is  an equivalence relation. We denote a cohomology classes of a continuous function $g$ by $[g]$. Then 
\begin{align}\nonumber
V:=\lbrace[g]:g:X\rightarrow\mathbb{R}\ \text{be\ continuous}\rbrace  
\end{align}
is vector space with respect to well-defined vector operation defined by 
\begin{align}\nonumber
   [g_1]+[g_2]:=[g_1+g_2],\ c[g_1]:=[cg_1]\ \text{for\ }g_1,g_2\ \text{be\ continuous\ and}\ c\in\mathbb{R}. 
\end{align}
We say that two continuous function $g_1,g_2:X\rightarrow \mathbb{R}$ are linearly independent as cohomology classes if $[g_1],[g_2]$ are linearly independent.

Recall that $M_n:=\bigcup_{i=1}^n\lbrace\gamma^ih:\gamma\in\Gamma_0,h\in H_0\rbrace\cup H_0$ and $K_n:=M_n^{\mathbb{N}}\cap\Sigma_A$
for $n\in\mathbb{N}$. We define the function $a_j|_{\pi(K_n)}:\pi(K_n)\rightarrow \mathbb{R}$ given by $a_j|_{\pi(K_n)}(x)=a_j(x)$ for all $x\in \pi(K_n) $ for all $j\in\{1,2,\cdots,m\}$, $n\in\mathbb{N}$ and $M_{\infty}:=\Lambda_c(G)$.

\begin{lemma}\label{cohology}
There exists $N\geq1$ such that $a_1|_{\pi(K_n)},a_2|_{\pi(K_n)},\cdots,a_m|_{\pi(K_n)}$ are linearly independent as cohomology classes for all $n\geq N$. In particular, $a_1,a_2,\cdots,a_m$ are linearly independent as cohomology classes.
\end{lemma}

\begin{pf}
For all $l\in\mathbb{N}$ and $i\in\{1,2,\cdots,m\}$ we define a point $\boldsymbol{x}_{l,i}=(x_{l,1},\cdots,x_{l,m})\in\mathbb{R}_{\geq 0}^m$ by $x_{l,i}=l/(2m)$ and $x_{l,k}=0$ for $k\in\{1,2,\cdots,m\}$ with $k\neq i$. We take a large natural number $N\in\mathbb{N}$ such that $(1/(2m+1),\cdots,1/(2m+1))\in\text{Conv}\{x_{N,i}:1\leq i\leq m\}$. For a contradiction, we assume that there exists $n\geq N$ such that
\begin{align}\nonumber
    c_1[a_1|_{\pi(K_n)}]+c_2[a_2|_{\pi(K_n)}]+\cdots+c_m[a_m|_{\pi(K_n)}]=[0]\ (c_1,\cdots,c_m\in\mathbb{R})
\end{align}
for some $(c_1,c_2,\cdots,c_m)\neq0$.
 Thus, there exists a bounded continuous function $u$ such that 
    $\sum_{i=1}^mc_ia_i\circ\pi|_{\pi(K_n)}=u-u\circ\tilde f$. Thus, for all $\mu\in M(\tilde f)$ we have
    \begin{align}\label{only}
    \int\sum_{i=1}^{m} c_ia_i|_{\pi(K_n)}d\mu=\int (u-u\circ\tilde f)d\mu =0.
    \end{align}
On the other hand, by proof of first part of Lemma \ref{positive approximation} and choice of $N$, there exists $\tilde f$-invariant measures $\mu_1,\mu_2$ supported on $\pi(K_n)$ such that 
\begin{align}
  &\int a_i|_{\pi(K_n)} d\mu_1< \frac{1}{2(m+1)}\ \text{for}\ i\in \lbrace 1,\cdots,m\rbrace\ \text{with}\ c_i>0\nonumber,\\
   &\int a_i|_{\pi(K_n)} d\mu_1>\frac{1}{2(m+1)}\ \text{for}\ i \in \lbrace 1,\cdots,m\rbrace\ \text{with}\ c_i<0\nonumber\nonumber,\\
  &\int a_i|_{\pi(K_n)} d\mu_2> \frac{1}{2(m+1)}\ \text{for}\ i \in \lbrace 1,\cdots,m\rbrace\ \text{with}\ c_i>0\nonumber,\\
  &\int a_i|_{\pi(K_n)} d\mu_2<\frac{1}{2(m+1)}\ \text{for}\ i \in \lbrace 1,\cdots,m\rbrace\ \text{with}\ c_i<0\nonumber.
\end{align}
Therefore, we get 
\begin{align*}
   \sum_{i=1}^m c_i\int a_i|_{\pi(K_n)} d\mu_1<\sum_{i=1}^m\frac{c_i}{2(m+1)}\ \ \text{and\ }\ 
   \sum_{i=1}^m c_i\int a_i|_{\pi(K_n)} d\mu_2>\sum_{i=1}^m\frac{c_i}{2(m+1)}.
\end{align*} 
This implies that $\sum_{i=1}^m{c_i}/{2(m+1)}<0<\sum_{i=1}^m{c_i}/{2(m+1)}$ by (\ref{only}). This is a contradiction.
$\qed$
\end{pf}

\begin{prop}\label{prop regularity of the dimension spectrum}
The dimension spectrum $\boldsymbol{\alpha} \mapsto b(\boldsymbol{\alpha})$ is real analytic on $\mathbb{R}_{>0}^m$.
\end{prop}

\begin{pf}
By Proposition \ref{prop more important proposition}, there exists $\boldsymbol{q}(\boldsymbol{\alpha})\in\mathbb{R}_{>0}^m$ such that 
\begin{align}\nonumber
    p(\boldsymbol{\alpha},\boldsymbol{q}(\boldsymbol{\alpha}),b(\boldsymbol{\alpha}))=0\ \text{and}\ \frac{\partial}{\partial q_i}p(\boldsymbol{\alpha},\boldsymbol{q}(\boldsymbol{\alpha}),b(\boldsymbol{\alpha}))=0\ (i=1,2,\cdots,m).\nonumber
\end{align}
Define the map 
$F:\mathbb{R}^m_{>0}\times\mathbb{R}^m_{>0}\times\mathbb{R}_{\geq 0}\rightarrow\mathbb{R}^{m+1}$ by
\begin{align}
F(\boldsymbol{\alpha},\boldsymbol{q},b):=\left(p(\boldsymbol{\alpha},\boldsymbol{q},b), \frac{\partial}{\partial q_1}p(\boldsymbol{\alpha},\boldsymbol{q},b),\cdots, \frac{\partial}{\partial q_m}p(\boldsymbol{\alpha},\boldsymbol{q},b)\right)\nonumber.
\end{align}
We want to apply the implicit function theorem. Since we have $F(\boldsymbol{\alpha},\boldsymbol{q}(\boldsymbol{\alpha}),b(\boldsymbol{\alpha}))=0$, it is sufficient to show that the matrix 
\begin{align}\nonumber
\begin{pmatrix}
\frac{\partial}{\partial b}p(\boldsymbol{\alpha},\boldsymbol{q},b)
&\frac{\partial^2}{\partial b\partial q_1}p(\boldsymbol{\alpha},\boldsymbol{q},b)
&\cdots
&\frac{\partial^2}{\partial b\partial q_m}p(\boldsymbol{\alpha},\boldsymbol{q},b)\\
\frac{\partial}{\partial q_1}p(\boldsymbol{\alpha},\boldsymbol{q},b)
&\frac{\partial^2}{\partial q_1\partial q_1}p(\boldsymbol{\alpha},\boldsymbol{q},b)
&\cdots
&\frac{\partial^2}{\partial q_1\partial q_m}p(\boldsymbol{\alpha},\boldsymbol{q},b)\\
\vdots&\vdots&\ddots&\vdots\\
\frac{\partial}{\partial q_m}p(\boldsymbol{\alpha},\boldsymbol{q},b)
&\frac{\partial^2}{\partial q_m\partial q_1}p(\boldsymbol{\alpha},\boldsymbol{q},b)
&\cdots
&\frac{\partial^2}{\partial q_m\partial q_m}p(\boldsymbol{\alpha},\boldsymbol{q},b)\\
\end{pmatrix}    
\end{align}
is invertible at $(\boldsymbol{\alpha},\boldsymbol{q}(\boldsymbol{\alpha}),b(\boldsymbol{\alpha}))$. We have 
\begin{align}\nonumber
    \frac{\partial}{\partial q_i}p(\boldsymbol{\alpha},\boldsymbol{q}(\boldsymbol{\alpha}),b(\boldsymbol{\alpha}))=0\ (i=1,2,\cdots,m)\ \text{and}\ \nonumber
    \frac{\partial}{\partial b}p(\boldsymbol{\alpha},\boldsymbol{q}(\boldsymbol{\alpha}),b(\boldsymbol{\alpha}))=\lambda(\mu_{\boldsymbol{q}(\boldsymbol{\alpha})})>0\nonumber
\end{align}
where $\mu_{\boldsymbol{q}(\boldsymbol{\alpha})}$ denotes the equilibrium measure for the potential $\langle\boldsymbol{q}(\boldsymbol{\alpha}),(\Phi+\boldsymbol{\alpha})\rangle-b(\boldsymbol{\alpha})\log|\tilde f'|$. Therefore, it is sufficient to prove that the matrix 
\begin{align}\nonumber H:=
\begin{pmatrix}
\frac{\partial^2}{\partial q_1\partial q_1}p(\boldsymbol{\alpha},\boldsymbol{q},b)
&\cdots
&\frac{\partial^2}{\partial q_1\partial q_m}p(\boldsymbol{\alpha},\boldsymbol{q},b)\\
\vdots&\ddots&\vdots\\
\frac{\partial^2}{\partial q_m\partial q_1}p(\boldsymbol{\alpha},\boldsymbol{q},b)
&\cdots
&\frac{\partial^2}{\partial q_m\partial q_m}p(\boldsymbol{\alpha},\boldsymbol{q},b)\\
\end{pmatrix}    
\end{align}
is invertible at $(\boldsymbol{\alpha},\boldsymbol{q}(\boldsymbol{\alpha}),b(\boldsymbol{\alpha}))$.
 By Lemma \ref{cohology}, the matrix $H$ is positive definite. Since $H$ is a symmetric matrix, this implies that $H$ is invertible. By implicit function theorem, the dimension spectrum $\boldsymbol{\alpha}\rightarrow b(\boldsymbol{\alpha})$ is real analytic on $\mathbb{R}_{>0}^m$.
$\qed$    
\end{pf}

Recall that the irregular set $J_{\text{ir}}$ is defined by 
\begin{align*}
    J_{\text{ir}}:=\left\{x\in\Lambda_c(G):\exists i\in\{1,2,\cdots m\}\ \text{s.t.\ }\lim_{n\to\infty}\frac{1}{n}\sum_{j=0}^{n-1}a_{i,j}\ \text{does\ not\ exist}\right\}.
\end{align*}
By Lemma \ref{cohology} and {\cite[Theorem 7.5]{barreira2000sets}}, we obtain the following lemma which is used to show that the irregular set has the same Hausdorff dimension as $\Lambda_c(G)$.
\begin{lemma}\label{irregular set}
    There exists $N\geq 1$ such that for all $n\geq N$ we have $\dim_H(\pi(K_n))=\dim_H(J_{\text{ir}}\cap \pi(K_n))$.
\end{lemma}

\begin{prop}
    We have $\dim_H(J_{\text{ir}})=\dim_H(\Lambda_c(G))$.
\end{prop}
\begin{pf}
    Since $J_{\text{ir}}\subset\Lambda_c(G)$, we have $\dim_H(J_{\text{ir}})\leq\dim_H(\Lambda_c(G))$. Thus, we shall show that $\dim_H(J_{\text{ir}})\geq\dim_H(\Lambda_c(G))$. By Proposition \ref{approximation property}, for all $s\in\mathbb{R}$ we have $P(-s\log|\tilde f'|)=\lim_{n\to\infty}P_{\pi(K_n)}(-s\log|\tilde f'|)$. By Lemma \ref{irregular set} and Bowen's formula, we obtain
    \begin{align*}
    P_{\pi(K_n)}(-\dim_H(J_{\text{ir}}\cap \pi(K_n))\log|\tilde f'|)=P_{\pi(K_n)}(-\dim_H(\pi(K_n))\log|\tilde f'|)=0.
    \end{align*}
    Since $J_{\text{ir}}\cap \pi(K_n)\subset J_{\text{ir}}$, we have $P_{\pi(K_n)}(-\dim_H(J_{\text{ir}})\log|\tilde f'|)\leq P_{\pi(K_n)}(-\dim_H(J_{\text{ir}}\cap \pi(K_n))\log|\tilde f'|)=0$. Therefore, we obtain $\lim_{n\to\infty}P_{\pi(K_n)}(-\dim_H(J_{\text{ir}})\log|\tilde f'|)=P(-\dim_H(J_{\text{ir}})\log|\tilde f'|)\leq0$. On the other hand, by Theorem \ref{useful Bowen}, we have $P(-\dim_H(\Lambda_c(G))\log|\tilde f'|)=0$. This means that $\dim_H(J_{\text{ir}})\geq\dim_H(\Lambda_c(G))$.  $\qed$
\end{pf}

\section{Proof of Proposition \ref{Proposition m=1 case}}\label{section proof 2} 
In this section we prove Proposition \ref{Proposition m=1 case}. Let $m=1$. 
\begin{pf}
Put $s:=\dim_H(\Lambda_c(G))$. Let $0<\alpha_1<\alpha_2<\infty$. By Theorem \ref{useful Bowen} we have $P(-s\log|\tilde f'|)=0$. Thus, by Proposition \ref{approximation property}, for all $l \in \mathbb{N}$ there exists $n(l)\in\mathbb{N}$ such that $0=P(-s\log|\tilde f'|)\geq P_{\pi(K_{n(l)})}(-s\log|\tilde f'|)>-{1}/{l}$. Put $P_{\pi(K_{n(l)})}(-s\log|\tilde f'|):=-\eta_{l}\ (\eta_{l}>0)$. Without loss of generality we may assume that
$(n(l)_{l\in\mathbb{N}}$ is a strictly increasing sequence of $\mathbb{N}$.
There exists equilibrium measure $\mu_{l}$ for $P_{\pi(K_{n(l)})}(-s\log|\tilde f' |)$. Thus, we obtain $h(\mu_{l})-s \lambda(\mu_{l})=-\eta_{l}$.
On the other hand, since the measures $\mu_l$ is a Gibbs measure, we have
\begin{align*}
\lambda(\mu_{l})=\int \log|\tilde f'| d \mu_l\asymp \sum_{\gamma \in \Gamma_0, h \in H_0}\sum_{j=1}^{n(l)}2\log j \cdot \exp(-P_{K_l}(-s\log|\tilde f'|)) \cdot j^{-2s}
\geq e^{-1} \sum_{\gamma \in \Gamma_0, h \in H_0}\sum_{j=1}^{n(l)}2\log j  \cdot j^{-2s}.
\end{align*}
  Hence, we obtain $\liminf_{l\to\infty}\lambda(\mu_l)>0$ and 
\begin{align}\label{eq limit 1}
    \lim_{l \to \infty}\frac{h(\mu_{l})}{\lambda(\mu_{l})}=\lim_{l \to \infty}\left(s - \frac{\eta_{l}}{\lambda(\mu_{l})}\right) =s
\end{align}
Moreover, we have 
\begin{align*}
 &\int a_1 d \mu_{l}
 =\sum_{\gamma \in \Gamma_0, h \in H_0}\sum_{j= 1}^{n(l)}\int_{[\gamma^jh]} a_1 d\mu_{l} 
 = \sum_{\gamma \in \Gamma_0, h \in H_0}\sum_{j= 1}^{n(l)} j \mu_{l}([\gamma^jh])\\
 &\asymp \sum_{\gamma \in \Gamma_0, h \in H_0}\sum_{j= 1}^{n(l)} j \cdot \exp(-P_{K_l}(-s\log|\tilde f'|))\cdot \exp(\sup_{[\gamma^jh]}-s\log |\tilde f'|)
 \asymp \sum_{\gamma \in \Gamma_0, h \in H_0}\sum_{j= 1}^{n(l)} j^{-2s+1}.
\end{align*}
Since $\lim_{l \to \infty}\sum_{\gamma \in \Gamma_0, h \in H_0}\sum_{j= 1}^{n(l)} j^{-2s+1} = \infty$, we obtain 
\begin{align}\label{eq limit 2}
\lim_{l \to \infty} \int a_1 d \mu_l=\infty.
\end{align}
The two limits (\ref{eq limit 1}) and (\ref{eq limit 2}) mean that $\lim_{\alpha\to\infty}\dim_H(J(\alpha))=1$. 
Also, by (\ref{eq limit 1}) and (\ref{eq limit 2}), there exists $\mu_1 \in M(\tilde f)$ such that $\lambda(\mu_1)<\infty$, $\int a_1 d \mu_1 > \alpha_2$, and ${h(\mu_1)}/{\lambda(\mu_1)}>b(\alpha_1)$. On the other hand, since $b(\alpha_1)=\max \left\{ {h(\mu)}/{\lambda(\mu)}: \mu \in  M(\tilde f), \lambda(\mu)<\infty ,\int a_1 d\mu=\alpha_1 \right\}$, there exists $\mu_2 \in M(\tilde f)$ such that ${h(\mu_2)}/{\lambda(\mu_2)}=b(\alpha_1)$ and $\int a_1 d\mu_2=\alpha_1$. Since $\alpha_1<\alpha_2<\int a_1 d\mu_1$, there exists $p\in(0,1)$ such that $\alpha_2 =p\alpha_1 +(1-p)\int a_1 d\mu_1$. Put $\nu:=p\mu_2+(1-p)\mu_1$. Then, we have $\int a_1 d\nu=\alpha_2$, $\lambda(\nu)<\infty$, and 
\begin{align}
{h(\nu)}=ph(\mu_2)+(1-p)h(\mu_1)\nonumber
&>pb(\alpha_1)\lambda(\mu_2)+(1-p)b(\alpha_1)\lambda(\mu_1) \nonumber
= b(\alpha_1)\lambda(\nu). \nonumber
\end{align}
Thus, we obtain $b(\alpha_2)\geq {h(\nu)}/{\lambda(\nu)}>b(\alpha_1)$ 
and $b(\alpha_2)> b(\alpha_1).$ $\qed$
\end{pf}

\textbf{Acknowledgments}
I would like to express my gratitude to my 
advisor, Johannes Jaerisch, for engaging in valuable discussions with me throughout the creation of this paper.
 \bibliographystyle{abbrv}
\bibliography{reference}

\begin{thebibliography}{10}

\bibitem{barany2021birkhoff}
B.~B{\'a}r{\'a}ny, T.~Jordan, A.~K{\"a}enm{\"a}ki, and M.~Rams.
\newblock Birkhoff and {L}yapunov spectra on planar self-affine sets.
\newblock {\em International Mathematics Research Notices},
  2021(10):7966--8005, 2021.

\bibitem{barreira2001variational}
L.~Barreira and B.~Saussol.
\newblock Variational principles and mixed multifractal spectra.
\newblock {\em Transactions of the American Mathematical Society},
  353(10):3919--3944, 2001.

\bibitem{barreira2002higher}
L.~Barreira, B.~Saussol, and J.~Schmeling.
\newblock Higher-dimensional multifractal analysis.
\newblock {\em Journal de math{\'e}matiques pures et appliqu{\'e}es},
  81(1):67--91, 2002.

\bibitem{barreira2000sets}
L.~Barreira and J.~Schmeling.
\newblock Sets of “non-typical” points have full topological entropy and
  full {H}ausdorff dimension.
\newblock {\em Israel Journal of Mathematics}, 116(1):29--70, 2000.

\bibitem{beardon2012geometry}
A.~F. Beardon.
\newblock {\em The geometry of discrete groups}, volume~91.
\newblock Springer Science \& Business Media, 2012.

\bibitem{borthwick2007spectral}
D.~Borthwick.
\newblock {\em Spectral theory of infinite-area hyperbolic surfaces}.
\newblock Springer, 2007.

\bibitem{bowen1979markov}
R.~Bowen and C.~Series.
\newblock Markov maps associated with {F}uchsian groups.
\newblock {\em Publications Math{\'e}matiques de l'IH{\'E}S}, 50:153--170,
  1979.

\bibitem{climenhaga2012topological}
V.~Climenhaga.
\newblock Topological pressure of simultaneous level sets.
\newblock {\em Nonlinearity}, 26(1):241, 2012.

\bibitem{Dal}
F.~Dal’Bo.
\newblock {\em Geodesic and horocyclic trajectories}.
\newblock Springer-Verlag London, Ltd.London, 2011.

\bibitem{falconer2004fractal}
K.~Falconer.
\newblock {\em Fractal geometry: mathematical foundations and applications}.
\newblock John Wiley \& Sons, 2004.

\bibitem{fan2015multifractal}
A.-H. Fan, T.~Jordan, L.~Liao, and M.~Rams.
\newblock Multifractal analysis for expanding interval maps with infinitely
  many branches.
\newblock {\em Transactions of the American Mathematical Society},
  367(3):1847--1870, 2015.

\bibitem{iommi2015multifractale}
G.~Iommi and T.~Jordan.
\newblock Multifractal analysis for quotients of {B}irkhoff sums for countable
  markov maps.
\newblock {\em International Mathematics Research Notices}, 2015(2):460--498,
  2015.

\bibitem{iommi2015multifractal}
G.~Iommi and T.~Jordan.
\newblock Multifractal analysis of {B}irkhoff averages for countable {M}arkov
  maps.
\newblock {\em Ergodic Theory and Dynamical Systems}, 35(8):2559--2586, 2015.

\bibitem{iommi2017transience}
G.~Iommi, T.~Jordan, and M.~Todd.
\newblock Transience and multifractal analysis.
\newblock {\em Annales de l'Institut Henri Poincar{\'e} C}, 34(2):407--421,
  2017.

\bibitem{Jaerisch2016AMA}
J.~Jaerisch, M.~Kessebohmer, and S.~Munday.
\newblock A multifractal analysis for cuspidal windings on hyperbolic surfaces.
\newblock {\em Stochastics and Dynamics}, 2016.

\bibitem{jaerisch2021mixed}
J.~Jaerisch and H.~Takahasi.
\newblock Mixed multifractal spectra of {B}irkhoff averages for non-uniformly
  expanding one-dimensional {M}arkov maps with countably many branches.
\newblock {\em Advances in Mathematics}, 385:107778, 2021.

\bibitem{jaerisch2022multifractal}
J.~Jaerisch and H.~Takahasi.
\newblock Multifractal analysis of homological growth rates for hyperbolic
  surfaces.
\newblock {\em arXiv preprint arXiv:2204.08907}, 2022.

\bibitem{jenkinson2005zero}
O.~Jenkinson, R.~D. Mauldin, and M.~Urba{\'n}ski.
\newblock Zero temperature limits of {G}s-equilibrium states for countable
  alphabet subshifts of finite type.
\newblock {\em Journal of Statistical Physics}, 119:765--776, 2005.

\bibitem{johansson2008multifractal}
A.~Johansson, T.~Jordan, A.~{\"O}berg, and M.~Pollicott.
\newblock Multifractal analysis of non-uniformly hyperbolic systems.
\newblock {\em Isr. J. Math. 177}, 2010.

\bibitem{Jordan2017BirkhoffSF}
T.~Jordan and M.~Rams.
\newblock Birkhoff spectrum for piecewise monotone interval maps.
\newblock {\em Fundamenta Mathematicae}, 2017.

\bibitem{kato2013perturbation}
T.~Kato.
\newblock {\em Perturbation theory for linear operators}, volume 132.
\newblock Springer Science \& Business Media, 2013.

\bibitem{katok1992fuchsian}
S.~Katok.
\newblock {\em Fuchsian groups}.
\newblock University of Chicago press, 1992.

\bibitem{Kessebhmer2007HigherdimensionalMV}
M.~Kesseb{\"o}hmer and M.~Urbanski.
\newblock Higher-dimensional multifractal value sets for conformal infinite
  graph directed markov systems.
\newblock {\em Nonlinearity}, 20:1969 -- 1985, 2007.

\bibitem{mauldin2003graph}
R.~D. Mauldin and M.~Urbanski.
\newblock {\em Graph directed {M}arkov systems: geometry and dynamics of limit
  sets}, volume 148.
\newblock Cambridge University Press, 2003.

\bibitem{Munday2011OnHD}
S.~Munday.
\newblock On {H}ausdorff dimension and cusp excursions for fuchsian groups.
\newblock {\em Discrete Contin. Dyn. Syst}, 2012.

\bibitem{przytycki2010conformal}
F.~Przytycki and M.~Urba{\'n}ski.
\newblock {\em Conformal fractals: ergodic theory methods}, volume 371.
\newblock Cambridge University Press, 2010.

\bibitem{ruelle2004thermodynamic}
D.~Ruelle.
\newblock {\em Thermodynamic formalism: the mathematical structure of
  equilibrium statistical mechanics}.
\newblock Cambridge University Press, 2004.

\bibitem{rush2023multifractal}
T.~Rush.
\newblock Multifractal analysis for {M}arkov interval maps with countably many
  branches.
\newblock {\em Nonlinearity}, 36(4):2038, 2023.

\bibitem{sarig2003existence}
O.~Sarig.
\newblock Existence of gibbs measures for countable {M}arkov shifts.
\newblock {\em Proceedings of the American Mathematical Society},
  131(6):1751--1758, 2003.

\bibitem{sarig2009lecture}
O.~Sarig.
\newblock Lecture notes on thermodynamic formalism for topological {M}arkov
  shifts.
\newblock {\em Penn State}, 2009.

\bibitem{sarig1999thermodynamic}
O.~M. Sarig.
\newblock Thermodynamic formalism for countable {M}arkov shifts.
\newblock {\em Ergodic Theory and Dynamical Systems}, 19(6):1565--1593, 1999.

\bibitem{series1985modular}
C.~Series.
\newblock The modular surface and continued fractions.
\newblock {\em Journal of the London Mathematical Society}, 2(1):69--80, 1985.

\bibitem{takahasi2020entropy}
H.~Takahasi.
\newblock Entropy-approachability for transitive markov shifts over infinite
  alphabet.
\newblock {\em Proceedings of the American Mathematical Society},
  148(9):3847--3857, 2020.

\bibitem{tanaka}
H.~Tanaka.
\newblock Higher-order asymptotic behaviours of pressure functionals and
  statistical representations of the coefficients.
\newblock {\em preprint}, 2022.

\bibitem{urubanskinoninvertible}
M.~Urba{\'n}ski, M.~Roy, and S.~Munday.
\newblock {\em Non-Invertible Dynamical Systems: Volume 2 Finer Thermodynamic
  Formalism–Distance Expanding Maps and Countable State Subshifts of Finite
  Type, Conformal GDMSs, Lasota-Yorke Maps and Fractal Geometry}, volume 490.
\newblock De Gruyter Expositions in Mathematics, 2022.

\bibitem{walters2000introduction}
P.~Walters.
\newblock {\em An introduction to ergodic theory}, volume~79.
\newblock Springer Science \& Business Media, 2000.

\end{thebibliography}
 \nocite{*}

\end{document}